\documentclass[12pt]{amsart}
\usepackage{amssymb,amsmath,amsfonts,latexsym,setspace,amscd}
\usepackage{bm}
\newcommand{\bea}{\begin{eqnarray}}
\newcommand{\eea}{\end{eqnarray}}

\newcommand{\cla}{\mathcal{A}}
\newcommand{\clb}{\mathcal{B}}

\newcommand{\cld}{\mathcal{D}}
\newcommand{\cle}{\mathcal{E}}
\newcommand{\clf}{\mathcal{F}}

\newcommand{\clh}{\mathcal{H}}
\newcommand{\clk}{\mathcal{K}}
\newcommand{\cll}{\mathcal{L}}
\newcommand{\clm}{\mathcal{M}}
\newcommand{\cln}{\mathcal{N}}
\newcommand{\clo}{\mathcal{O}}

\newcommand{\clq}{\mathcal{Q}}
\newcommand{\clr}{\mathcal{R}}
\newcommand{\cls}{\mathcal{S}}

\newcommand{\clu}{\mathcal{U}}
\newcommand{\clw}{\mathcal{W}}

\newcommand{\clz}{\mathcal{Z}}
\newcommand{\z}{\bm{z}}
\newcommand{\w}{\bm{w}}

\newcommand{\raro}{\rightarrow}

\def \qed {\hfill \vrule height6pt width 6pt depth 0pt}
\def\textmatrix#1&#2\\#3&#4\\{\bigl({#1 \atop #3}\ {#2 \atop #4}\bigr)}
\def\dispmatrix#1&#2\\#3&#4\\{\left({#1 \atop #3}\ {#2 \atop #4}\right)}
\newcommand{\be}{\begin{equation}}
\newcommand{\ee}{\end{equation}}
\newcommand{\ben}{\begin{eqnarray*}}
\newcommand{\een}{\end{eqnarray*}}

\newcommand{\NI}{\noindent}

\newcommand{\bi}{\begin{itemize}}
\newcommand{\ei}{\end{itemize}}

\newtheorem{Theorem}{\sc Theorem}[section]
\newtheorem{Lemma}[Theorem]{\sc Lemma}
\newtheorem{Proposition}[Theorem]{\sc Proposition}
\newtheorem{Corollary}[Theorem]{\sc Corollary}
\newtheorem{Definition}[Theorem]{\sc Definition}
\newtheorem{Example}[Theorem]{\sc Example}
\newtheorem{Remark}[Theorem]{\sc Remark}
\newtheorem{Note}[Theorem]{\sc Note}
\newtheorem{Question}{\sc Question}
\newtheorem{ass}[Theorem]{\sc Assumption}
\newcommand{\bt}{\begin{Theorem}}
\def\beginlem{\begin{Lemma}}
\def\beginprop{\begin{Proposition}}
\def\begincor{\begin{Corollary}}
\def\begindef{\begin{Definition}}
\def\beginexamp{\begin{Example}}
\def\beginrem{\begin{Remark}}
\def\beginq{\begin{Question}}
\def\beginass{\begin{ass}}
\def\beginnote{\begin{Note}}
\newcommand{\et}{\end{Theorem}}
\def\endlem{\end{Lemma}}
\def\endprop{\end{Proposition}}
\def\endcor{\end{Corollary}}
\def\enddef{\end{Definition}}
\def\endexamp{\end{Example}}
\def\endrem{\end{Remark}}
\def\endq{\end{Question}}
\def\endass{\end{ass}}
\def\endnote{\end{Note}}

\begin{document}

\title[Hilbert Module Approach to Multivariable Operator Theory]{An Introduction to Hilbert Module Approach to Multivariable Operator Theory}

\author[Jaydeb Sarkar]{Jaydeb Sarkar}
\address{Indian Statistical Institute, Statistics and Mathematics Unit, 8th Mile, Mysore Road, Bangalore, 560059, India}
\email{jay@isibang.ac.in, jaydeb@gmail.com}




\begin{abstract}
Let $\{T_1, \ldots, T_n\}$ be a set of $n$ commuting bounded linear
operators on a Hilbert space $\clh$. Then the $n$-tuple $(T_1,
\ldots, T_n)$ turns $\clh$ into a module over $\mathbb{C}[z_1,
\ldots, z_n]$ in the following sense:
\[\mathbb{C}[z_1, \ldots, z_n] \times \clh \raro \clh, \quad \quad (p, h) \mapsto
p(T_1, \ldots, T_n)h,\]where $p \in \mathbb{C}[z_1, \ldots, z_n]$
and $h \in \clh$. The above module is usually called the Hilbert
module over $\mathbb{C}[z_1, \ldots, z_n]$. Hilbert modules over
$\mathbb{C}[z_1, \ldots, z_n]$ (or natural function algebras) were
first introduced by R. G. Douglas and C. Foias in 1976. The two main
driving forces were the algebraic and complex geometric views to
multivariable operator theory.

\NI This article gives an introduction of Hilbert modules over
function algebras and surveys some recent developments. Here the
theory of Hilbert modules is presented as combination of commutative
algebra, complex geometry and the geometry of Hilbert spaces and its
applications to the theory of $n$-tuples ($n \geq 1$) of commuting
operators. The topics which are studied include: model theory from
Hilbert module point of view, Hilbert modules of holomorphic
functions, module tensor products, localizations, dilations,
submodules and quotient modules, free resolutions, curvature and
Fredholm Hilbert modules. More developments in the study of Hilbert
module approach to operator theory can be found in a companion
paper, ``Applications of Hilbert Module Approach to Multivariable
Operator Theory''.
\end{abstract}

\maketitle

\tableofcontents


\section{Introduction}\label{I}

One of the most important areas of investigation in operator theory
is the study of $n$-tuples of commuting bounded linear operators on
Hilbert spaces, or Hilbert modules over natural function algebras. A
Hilbert module $\clh$ over $\mathbb{C}[z_1, \ldots, z_n]$ is the
Hilbert space $\clh$ equipped with $n$ module maps, that is, with an
$n$-tuple of commuting bounded linear operators on $\clh$.

The origins of Hilbert modules, in fact, lie in classical linear
operators on finite dimensional vector spaces. For instance, let $T$
be a linear operator on an $n$-dimensional vector space $\clh$. Then
$\clh$ is a module over $\mathbb{C}[z]$ in the following sense:
\[\mathbb{C}[z] \times \clh \raro \clh, \quad \quad (p, h) \mapsto
p(T)h, \quad \quad (p \in \mathbb{C}[z], h \in \clh)\]where for $p =
\sum_{k \geq 0} a_k z^k \in \mathbb{C}[z]$, $p(T)$ is the natural
functional calculus given by $p(T) = \sum_{k \geq 0} a_k T^k$. Since
$\mathbb{C}[z]$ is a principle ideal domain, the ideal $\{p \in
\mathbb{C}[z]: p(T) = 0\}$ is generated by a \textsf{non-zero}
polynomial. Such a polynomial is called a \textit{minimal
polynomial} for $T$. The existence of a minimal polynomial is a key
step in the classification of linear operators in finite dimensional
Hilbert spaces \cite{PF}. More precisely, the existence of the
Jordan form follows from the structure theorem for
finitely-generated modules over principle ideal domains.

The idea of viewing a commuting tuple of operators on a Hilbert
space as Hilbert module over a natural function algebra goes back to
Ronald G. Douglas in the middle of 1970. Perhaps, the main
motivations behind the Douglas approach to Hilbert modules were the
elucidating role of Brown-Douglas-Fillmore theory (1973), complex
geometric interpretation of the Cowen-Douglas class (1978),
Hormandar's algebraic approach, in the sense of Koszul complex, to
corona problem (1967) and later, Taylor's notion of joint spectrum
(1970), again in the sense of Koszul complex, in operator theory and
function theory.

Historically, the first ideas leading to Hilbert modules can be
traced back to the unpublished manuscript \cite{DF-76}, in which
Douglas and Foias proposed an algebraic approach to dilation theory.
Then in \cite{D1} and \cite{D2}, the notion of Hilbert modules
became refined.

A systematic study of Hilbert modules only really started in 1989
with the work of Douglas, Paulsen and Yan \cite{DPY} and the
monograph by Douglas and Paulsen \cite{DP}.

\NI Since then, this approach has become one of the essential tools
of multivariable operator theory. This field now has profound
connections to various areas of mathematics including commutative
algebra, complex geometry and topology (see \cite{D-KY},
\cite{DPSY}, \cite{KY}, \cite{DKKS2}, \cite{DKKS1}, \cite{DS1},
\cite{DS4}, \cite{DMS}).

The purpose of this survey article is to present a (Hilbert) module
approach to multivariable operator theory. The topics and results
covered here are chosen to complement the existing monographs by
Douglas and Paulsen \cite{DP} and Chen and Guo \cite{C-G} and
surveys by Douglas \cite{D-sur1} and \cite{D-sur2} though some
overlap will be unavoidable. Many interesting results, open problems
and references can be found in the monograph and the surveys
mentioned above.

In view of time and space constraints, the present survey will not
cover many interesting aspects of Hilbert module approach to
multivariable operator theory, including the case of single
operators. A few of these are: (1) The classification program for
reducing subspaces of multiplication by Blaschke products on the
Bergman space, by Zhu, Guo, Douglas, Sun, Wang, Putinar. (see
\cite{D10}, \cite{D11}, \cite{GH11}). (2) Extensions of Hilbert
modules by Carlson, Clark, Foias, Guo, Didas and Eschmeier (see
\cite{DE06}, \cite{CC1}, \cite{CC2}, \cite{G99}). (3) $K_0$-group
and similarity classification by Jiang, Wang, Ji, Guo (see
\cite{JJ}, \cite{JW}, \cite{Jia1}. (4) Classification programme of
homogeneous operators by Clark, Bagchi, Misra, Sastry and Koranyi
(see \cite{GM-NSN}, \cite{GM2}, \cite{GM3} and \cite{GM4}). (4)
Sheaf-theoretic techniques by Eschmeier, Albrecht, Putinar, Taylor
and Vasilescu (see \cite{EP}, \cite{Vas}).

Finally, although the main guiding principle of this development is
the correspondence
\[\mbox{commuting} \,n\mbox{-tuples} \quad
\longleftrightarrow \quad \mbox{Hilbert modules over~}
\mathbb{C}[z_1, \ldots, z_n],\]it is believed that the Hilbert
module approach is a natural way to understand the subject of
multivariable operator theory.

\NI\textsf{Outline of the paper:} The paper has seven sections
besides this introduction. Section 2 begins with a brief
introduction of Hardy module which is a well-established procedure
to pass from the function theory to the one variable operator
theory. This section also includes basics of Hilbert modules over
function algebras, localizations and dilations. The third section is
centered around those aspects of operator theory that played an
important role in the development of Hilbert modules. In particular,
the third section introduce three basic notions which are directly
formulated with the required structures, namely, algebraic, analytic
and geometric. Section 4 is devoted to the study of contractive
Hilbert modules over $\mathbb{C}[z]$. Section 5 describes the
relationship of von Neumann-Wold decomposition with the structure of
submodules of the Hardy (general function Hilbert) module(s).
Section 6 introduces the notion of unitarily equivalence submodules
of Hilbert modules of holomorphic functions. Section 7 sets the
homological framework for Hilbert modules and Section 8 introduces
the theory of Drury-Arveson module.

\NI\textsf{Notations and Conventions:} (i) $\mathbb{N} =$ Set of all
natural numbers including 0.  (ii) $n \in \mathbb{N}$ and $n \geq
1$, unless specifically stated otherwise. (iii) $\mathbb{N}^n =
\{\bm{k} = (k_1, \ldots, k_n) : k_i \in \mathbb{N}, i = 1, \ldots,
n\}$. (iv)  $\mathbb{C}^n =$ the complex $n$-space. (v) $\Omega$ :
Bounded domain in $\mathbb{C}^n$. (vi) $\bm{z} = (z_1, \ldots, z_n)
\in \mathbb{C}^n$. (vii) ${z}^{\bm{k}} = z_1^{k_1}\cdots z_n^{k_n}$.
(viii) $\clh, \clk, \cle, \cle_*$ : Hilbert spaces. (ix) $\clb(\clh,
\clk)=$ the set of all bounded linear operators from $\clh$ to
$\clk$. (x) $T = (T_1, \ldots, T_n)$, $n$-tuple of commuting
operators. (xi) $T^{\bm{k}} = T_1^{k_1} \cdots T_n^{k_n}$. (xii)
$\mathbb{C}[\z] = \mathbb{C}[z_1, \ldots, z_n]$. (xiii)
$\mathbb{D}^n = \{\z : |z_i| <1, i = 1, \ldots, n\}$, $\mathbb{B}^n
= \{\z : \|\z\|_{\mathbb{C}^n} <1\}$. (xiv) $H^2_{\cle}(\mathbb{D})$
: $\cle$-valued Hardy space over $\mathbb{D}$.

\NI Throughout this note all Hilbert spaces are over the complex
field and separable. Also for a closed subspace $\cls$ of a Hilbert
space $\clh$, the orthogonal projection of $\clh$ onto $\cls$ will
be denoted by $P_{\cls}$.

\section{Hilbert modules}

The purpose of this section is to give some of the essential
background for Hilbert modules. The first subsection is devoted to
set up the notion of Hilbert modules over the polynomial algebra.
The third subsection deals with Hilbert modules over function
algebras. Basic concepts and classical definitions are summarized in
the subsequent subsections.

Before proceeding to the detailed development, it is more convenient
to introduce a brief overview of the Hardy space over the unit disc
$\mathbb{D}$. Results based on the Hardy space and the
multiplication operator on the Hardy space play an important role in
both operator theory and function theory. More precisely, for many
aspects of geometric and analytic intuition, the Hardy space
techniques play a fundamental role in formulating problems in
operator theory and function theory both in one and several
variables.

The \textit{Hardy space} $H^2(\mathbb{D})$ over $\mathbb{D}$ is the
set of all power series \[f = \sum_{m=0}^{\infty} a_m z^m, \quad
\quad (a_m \in \mathbb{C})\]such that
\[\|f\|_{H^2(\mathbb{D})} := (\sum_{m=0}^\infty
|a_m|^2)^{\frac{1}{2}} < \infty.\]Let $f = \sum_{m=0}^{\infty} a_m
z^m \in H^2(\mathbb{D})$. It is obvious that $\sum_{m=0}^\infty
|w|^m < \infty$ for each $w \in \mathbb{D}$. This and
$\sum_{m=0}^\infty |a_m|^2 < \infty$ readily implies that
\[\sum_{m=0}^\infty a_m w^m\]converges absolutely for each $w \in \mathbb{D}$. In other words,
$f = \sum_{m=0}^{\infty} a_m z^m$ is in $H^2(\mathbb{D})$ if and
only if $f$ is a square summable holomorphic function on
$\mathbb{D}$.

\NI Now, for each $w \in \mathbb{D}$ one can define a complex-valued
function $\mathbb{S}(\cdot, w) : \mathbb{D} \raro \mathbb{C}$ by
\[(\mathbb{S}(\cdot, w))(z) = \sum_{m=0}^\infty \bar{w}^m z^m. \quad \quad (z \in
\mathbb{D})\]Since \[\sum_{m=0}^\infty |\bar{w}^m|^2 =
\sum_{m=0}^\infty (|{w}|^2)^m = \frac{1}{1 - |w|^2},\]it follows
that
\[\mathbb{S}(\cdot, w) \in H^2(\mathbb{D}),\quad \quad (w \in
\mathbb{D})\]and \[\|\mathbb{S}(\cdot, w)\|_{H^2(\mathbb{D})} =
\frac{1}{(1 - |w|^2)^{\frac{1}{2}}}. \quad \quad (w \in
\mathbb{D})\]Moreover, if $f = \sum_{m=0}^\infty a_m z^m \in
H^2(\mathbb{D})$ and $w \in \mathbb{D}$, then \[f(w) =
\sum_{m=0}^\infty a_m w^m = \langle \sum_{m=0}^\infty a_m z^m,
\sum_{m=0}^\infty \bar{w}^m z^m \rangle_{H^2(\mathbb{D})} = \langle
f, \mathbb{S}(\cdot, w)\rangle_{H^2(\mathbb{D})}.\]Therefore, the
vector $\mathbb{S}(\cdot, w) \in H^2(\mathbb{D})$
\textit{reproduces} (cf. Subsection \ref{sub-RKHM}) the value of $f
\in H^2(\mathbb{D})$ at $w \in \mathbb{D}$. In particular,
\[(\mathbb{S}(\cdot, w)) (z) = \langle \mathbb{S}(\cdot, w),
\mathbb{S}(\cdot, z)\rangle_{H^2(\mathbb{D})} = \sum_{m=0}^\infty
z^m \bar{w}^m = (1 - z \bar{w})^{-1}. \quad \quad (z, w \in
\mathbb{D})\]The function $\mathbb{S} : \mathbb{D} \times \mathbb{D}
\raro \mathbb{C}$ defined by \[\mathbb{S}(z, w) = (1 - z
\bar{w})^{-1}, \quad \quad (z, w \in \mathbb{D})\]is called the
\textit{Szeg\H{o}} or \textit{Cauchy-Szeg\H{o}} kernel of
$\mathbb{D}$. Consequently, $H^2(\mathbb{D})$ is a
\textit{reproducing kernel Hilbert space} with kernel function
$\mathbb{S}$ (see Subsection \ref{sub-RKHM}).

The next goal is to show that the set $\{\mathbb{S}(\cdot, w) : w
\in \mathbb{D}\}$ is \textit{total} in $H^2(\mathbb{D})$, that is,
\[\overline{\mbox{span}} \{\mathbb{S}(\cdot, w) : w \in \mathbb{D}\}
= H^2(\mathbb{D}).\]To see this notice that the reproducing property
of the Szeg\H{o} kernel yields $f(w) = \langle f, \mathbb{S}(\cdot,
w)\rangle_{H^2(\mathbb{D})}$ for all $f \in H^2(\mathbb{D})$ and $w
\in \mathbb{D}$. Now the result follows from the fact that \[f \perp
\mathbb{S}(\cdot, w),\]for $f \in H^2(\mathbb{D})$ and for all $w
\in \mathbb{D}$ if and only if
\[f = 0.\]It also follows that for each $w \in \mathbb{D}$, the
\textit{evaluation map} $ev_w : H^2(\mathbb{D}) \raro \mathbb{C}$
defined by
\[ev_w(f) = f(w), \quad \quad (f \in H^2(\mathbb{D}))\]is
continuous.

The next task is to recall some of the most elementary properties of
the multiplication operator on $H^2(\mathbb{D})$. Observe first that
\[\langle z(z^k), z(z^l)\rangle_{H^2(\mathbb{D})} = \langle z^{k+1},
z^{l+1}\rangle_{H^2(\mathbb{D})} = \delta_{k,l} = \langle z^k,
z^l\rangle_{H^2(\mathbb{D})}. \quad \quad (k, l \in \mathbb{N})\]
Using the fact that the set $\{z^m : m \in \mathbb{N}\}$ is total in
$H^2(\mathbb{D})$, the previous equality implies that the
multiplication operator $M_z$ on $H^2(\mathbb{D})$ defined by
\[(M_z f)(w) = w f(w),\quad \quad (f \in H^2(\mathbb{D}), w \in
\mathbb{D})\]is an isometric operator, that is, \[M_z^* M_z =
I_{H^2(\mathbb{D})}.\]Moreover, \[\langle M_z^* z^k,
z^l\rangle_{H^2(\mathbb{D})} = \langle z^k,
z^{l+1}\rangle_{H^2(\mathbb{D})} = \delta_{k, l+1} = \delta_{k-1, l}
= \langle z^{k-1}, z^l\rangle_{H^2(\mathbb{D})},\]for all $k \geq 1$
and $l \in \mathbb{N}$. Also it follows that $\langle M_z^* 1,
z^l\rangle_{H^2(\mathbb{D})} = 0$. Consequently,  \[ M_z^* z^k =
\left\{
\begin{array}{ll}
z^{k-1} & \;\;\mbox{if $k \geq 1$};\\
0 & \;\;\mbox{if $k = 0$}.\end{array} \right.\]It also follows that
\[\begin{split}\langle(I_{H^2(\mathbb{D})} - M_z M_z^*) \mathbb{S}(\cdot, w),
\mathbb{S}(\cdot, z)\rangle_{H^2(\mathbb{D})} & = \langle
\mathbb{S}(\cdot, w), \mathbb{S}(\cdot, z)\rangle_{H^2(\mathbb{D})}
- \langle M_z^* \mathbb{S}(\cdot, w), M_z^* \mathbb{S}(\cdot,
z)\rangle_{H^2(\mathbb{D})} \\ & = \mathbb{S}(z, w) - z \bar{w}
\mathbb{S}(z, w) = 1 \\ & = \langle P_{\mathbb{C}} \mathbb{S}(\cdot,
w), \mathbb{S}(\cdot, z)\rangle_{H^2(\mathbb{D})},\end{split}\]where
$P_{\mathbb{C}}$ is the orthogonal projection of $H^2(\mathbb{D})$
onto the one-dimensional subspace of all constant functions on
$\mathbb{D}$. Therefore, \[I_{H^2(\mathbb{D})} - M_z M_z^* =
P_{\mathbb{C}}.\]To compute the kernel, $\mbox{ker}(M_z - w
I_{H^2(\mathbb{D})})^*$ for $w \in \mathbb{D}$, note that
\[\begin{split}M_z^* \mathbb{S}(\cdot, w) & = M_z^*(1 + \bar{w} z + \bar{w}^2 z^2 + \cdots) = \bar{w} + \bar{w}^2 z + \bar{w}^3 z^2 + \cdots
= \bar{w} (1 + \bar{w} z + \bar{w}^2 z^2 + \cdots) \\ &= \bar{w}
\mathbb{S}(\cdot, w).\end{split}\] On the other hand, if $M_z^* f =
\bar{w} f$ for some $f \in H^2(\mathbb{D})$ then
\[f(0) = P_{\mathbb{C}} f = (I_{H^2(\mathbb{D})} - M_z M_z^*)f = (1 - z \bar{w})
f,\]that is, $f = f(0) \mathbb{S}(\cdot, w)$. Consequently, $M_z^* f
= \bar{w} f$ if and only if $f = \lambda \mathbb{S}(\cdot, w)$ for
some $\lambda \in \mathbb{C}$. That is,
\[\mbox{ker}(M_z - w I_{H^2(\mathbb{D})})^* = \{\lambda
\mathbb{S}(\cdot, w) : \lambda \in \mathbb{C}\}.\]In particular,
\[\mathop{\bigvee}_{w \in \mathbb{D}} \mbox{ker}(M_z - w
I_{H^2(\mathbb{D})})^* = H^2(\mathbb{D}).\]

The following theorem summarizes the above observations.

\begin{Theorem}
Let $H^2(\mathbb{D})$ denote the Hardy space over $\mathbb{D}$ and
$M_z$ denote the multiplication operator by the coordinate function
$z$ on $H^2(\mathbb{D})$. Then, the following properties hold:

(i) The set $\{\mathbb{S}(\cdot, w) : w \in \mathbb{D}\}$ is total
in $H^2(\mathbb{D})$.

(ii) The evaluation map $ev_w : H^2(\mathbb{D}) \raro \mathbb{C}$
defined by $ev_w(f) = f(w)$ is continuous for each $w \in
\mathbb{D}$.

(iii) $\sigma_p(M_z^*) = \mathbb{D}$ and $\mbox{ker}(M_z - w
I_{H^2(\mathbb{D})})^* = \{ \lambda \mathbb{S}(\cdot, w): \lambda
\in \mathbb{C}\}$.

(iv) $f(w) = \langle f, \mathbb{S}(\cdot,
w)\rangle_{H^2(\mathbb{D})}$ for all $f \in H^2(\mathbb{D})$ and $w
\in \mathbb{D}$.

(v)$I_{H^2(\mathbb{D})} - M_z M_z^* = P_{\mathbb{C}}$.

(vi) $\mathop{\bigvee}_{w \in \mathbb{D}} \mbox{ker}(M_z - w
I_{H^2(\mathbb{D})})^* = H^2(\mathbb{D})$.
\end{Theorem}

Let $\cle$ be a Hilbert space. In what follows,
$H^2_{\cle}(\mathbb{D})$ stands for the Hardy space of $\cle$-valued
analytic functions on $\mathbb{D}$. Moreover, by virtue of the
unitary $U : H^2_{\cle}(\mathbb{D}) \raro H^2(\mathbb{D}) \otimes
\cle$ defined by \[z^m \eta \mapsto z^m \otimes \eta, \quad (\eta
\in \cle, m \in \mathbb{N})\] the vector valued Hardy space
$H^2_{\cle}(\mathbb{D})$ will be identified with the Hilbert space
tensor product $H^2(\mathbb{D}) \otimes \cle$.

For a more extensive treatment of the Hardy space and related
topics, the reader is referred to the books by Sz.-Nagy and Foias
\cite{NF}, Radjavi and Rosenthal \cite{RR}, Rosenblum and Rovnyak
\cite{RoRo} and Halmos \cite{H-book}.

\subsection{Hilbert modules over $\mathbb{C}[\bm{z}]$}\label{HMC}

Let $\{T_1, \ldots, T_n\}$ be a set of $n$ commuting bounded linear
operators on a Hilbert space $\clh$. Then the $n$-tuple $(T_1,
\ldots, T_n)$ turns $\clh$ into a module over $\mathbb{C}[\bm{z}]$
in the following sense:
\[\mathbb{C}[\bm{z}] \times \clh \raro \clh, \quad \quad (p, h)
\mapsto p(T_1, \ldots, T_n)h,\]where $p \in \mathbb{C}[\bm{z}]$ and
$h \in \clh$. The above module is usually called the \textit{Hilbert
module} over $\mathbb{C}[\bm{z}]$. Denote by $M_p: \clh \raro \clh$
the bounded linear operator
\[
M_p h = p \cdot h = p(T_1, \ldots, T_n)h, \quad \quad(h \in \clh)
\]
for $p \in \mathbb{C}[\z]$. In particular, for $p = z_i \in
\mathbb{C}[\z]$, this gives the \textit{module multiplication}
operators $\{M_{j}\}_{j=1}^n$ by the coordinate functions
$\{z_j\}_{j=1}^n$ defined by
\[M_{i} h = z_i(T_1, \ldots, T_n) h = T_i h. \quad \quad(h \in
\clh,\, 1 \leq i \leq n)\] Here and in what follows, the notion of a
Hilbert module $\clh$ over $\mathbb{C}[\z]$ will be used in place of
an $n$-tuple of commuting operators $\{T_1, \ldots, T_n\} \subseteq
\clb(\clh)$, where the operators are determined by module
multiplication by the coordinate functions, and vice versa.

\textit{When necessary, the notation $\{M_{\clh,i}\}_{i=1}^n$ will
be used to indicate the underlying Hilbert space $\clh$ with respect
to which the module maps are defined.}

Let $\cls$ be a closed subspace of $\clh$. Then $\cls$ is a
\textit{submodule} of $\clh$ if $M_{i} \cls \subseteq \cls$ for all
$i = 1, \ldots, n$. A closed subspace $\clq$ of $\clh$ is said to be
\textit{quotient module} of $\clh$ if $\clq^\perp \cong \clh/ \clq$
is a submodule of $\clh$. Therefore, a closed subspace $\clq$ is a
quotient module of $\clh$ if and only if $M^*_{i} \clq \subseteq
\clq$ for all $i = 1, \ldots, n$. In particular, if the module
multiplication operators on a Hilbert module $\clh$ are given by the
commuting tuple of operators $(T_1, \ldots, T_n)$ then $\cls$ is a
submodule of $\clh$ if and only if $\cls$ is joint $(T_1, \ldots,
T_n)$-invariant subspace of $\clh$ and $\clq$ is a quotient module
of $\clh$ if and only if $\clq$ is joint $(T_1^*, \ldots,
T_n^*)$-invariant subspace of $\clh$.

Let $\cls$ be a submodule and $\clq$ be a quotient module of a
Hilbert module $\clh$ over $\mathbb{C}[\z]$. Then $\cls$ and $\clq$
are also Hilbert modules over $\mathbb{C}[\bm{z}]$ where the module
multiplication by the coordinate functions on $\cls$ and $\clq$ are
given by the restrictions $(R_{1}, \ldots, R_{n})$ and the
compressions $(C_{1}, \ldots, C_{n})$ of the module multiplication
operators on $\clh$, respectively. That is,
\[R_{i} = M_{i}|_{\cls} \quad \mbox{and} \quad C_{i} = P_{\clq}
M_{i}|_{\clq}. \quad \quad (1 \leq i \leq n)\]Evidently,
\[R_{i}^* = P_{\cls} M_{i}^*|_{\cls} \quad \quad \mbox{and}
\quad \quad C_{i}^* = M_{i}^*|_{\clq}. \quad \quad (1 \leq i \leq
n)\]

A bounded linear map $X: \clh \rightarrow \clk$ between two Hilbert
modules $\clh$ and $\clk$ over $\mathbb{C}[\bm{z}]$ is said to be a
\textit{module map} if $X M_{i} = M_{i} X$ for $i = 1, \ldots, n$,
or equivalently, if $X M_p = M_p X$ for $p \in \mathbb{C}[\bm{z}]$.
A pair of Hilbert modules will be considered the 'same', that is,
\textit{isomorphic} provided there is a unitary module map between
them, and \textit{similar} if there is an invertible module map
between them.

\subsection{Dilations} The purpose of this subsection is to present a modified version
of dilation theory for commuting tuples of operators.

Let $\clh$ and $\clk$ be Hilbert modules over $\mathbb{C}[\bm{z}]$.
Then

(1) A map $\Pi \in \clb(\clh, \clk)$ is called \textit{co-module}
map if $\Pi^* : \clk \raro \clh$ is a module map, that is, $\Pi
M_i^* = M_i^* \Pi$.

(2) $\clk$ is said to be \textit{dilation} of $\clh$ if there exists
a co-module isometry $\Pi : \clh \raro \clk$.  In this case, we also
say that $\Pi \in \clb(\clh, \clk)$ is a dilation of $\clh$.

(3) A dilation $\Pi \in \clb(\clh, \clk)$ of $\clh$ is
\textit{minimal} if $\clk = \overline{\mbox{span}} \{M^{\bm{k}} (\Pi
\clh) : \bm{k} \in \mathbb{N}^n\}$.

Let $\Pi \in \clb(\clh, \clk)$ be a dilation of $\clh$. Then
$\Pi(\clh)$ is a quotient module of $\clk$, that is, $\Pi(\clh)$ is
a joint $(M_{1}^*, \ldots, M_{n}^*)$-invariant subspace of $\clk$,
and
\[M^{\bm{k}} = P_{\Pi(\clh)} M^{\bm{k}}|_{\Pi(\clh)},\]for
all $\bm{k} \in \mathbb{N}^n$. Moreover, one has the following short
exact sequence of Hilbert modules
\begin{equation}\label{ses}0 \longrightarrow \cls \stackrel{i}
\longrightarrow \clk \stackrel{\pi} \longrightarrow \clh
\longrightarrow 0,\end{equation} where $\cls = (\Pi \clh)^{\perp}
(\cong \clk/ \Pi \clh)$, a submodule of $\clk$, $i$ is the inclusion
and $\pi:= \Pi^*$ is the quotient map. In other words, if $\clk$ is
a dilation of $\clh$ then there exists a quotient module $\clq$ and
a submodule $\cls$ of $\clk$ such that $\clk = \cls \oplus \clq$,
that is, \[0 \longrightarrow \cls \stackrel{i} \longrightarrow \clk
\stackrel{\pi} \longrightarrow \clq \longrightarrow 0,\]and $\clq
\cong \clh$.

\NI Conversely, let $\clh$ and $\clk$ be Hilbert modules over
$\mathbb{C}[\z]$ and $\clh \cong \clq$, a quotient module of $\clk$.
Therefore, $\clk$ is a dilation of $\clh$ and by defining $\cls :=
\clq^{\perp}$, a submodule of $\clk$, one arrives at the short exact
sequence (\ref{ses}).

A Hilbert module $\clh$ over $\mathbb{C}[z]$ is said to be
\textit{contractive Hilbert module over $\mathbb{C}[z]$} if
$I_{\clh} - M^* M \geq 0$.

The famous isometric dilation theorem of Sz. Nagy (cf. \cite{NF})
states that:

\begin{Theorem}\label{Nagy-dilation}
Every contractive Hilbert module $\clh$ over $\mathbb{C}[z]$ has a
minimal isometric dilation.
\end{Theorem}

\NI \textsf{Proof.} Let $\clh$ be a contractive Hilbert module over
$\mathbb{C}[z]$. Let $D_{\clh} = (I_{\clh} - M^* M)^{\frac{1}{2}}$
and $\cln_{\clh} := \clh \oplus H^2_{\clh}(\mathbb{D})$. Define
$N_{\clh} \in \clb(\cln_{\clh})$ by
\[N_{\clh} :=
\begin{bmatrix}M& 0\\\bm{D}_{\clh}&M_z\end{bmatrix},\]where
$\bm{D}_{\clh} : \clh \raro H^2_{\clh}(\mathbb{D})$ is the constant
function defined by $(\bm{D}_{\clh} h)(z) = D_{\clh}h$ for all $h
\in \clh$ and $z \in \mathbb{D}$. Consequently, \[N_{\clh}^*
N_{\clh} = \begin{bmatrix}M^* M + D^2_{\clh}&
0\\0&I_{H^2_{\clh}(\mathbb{D})}\end{bmatrix} =
\begin{bmatrix}I_{\clh}&
0\\0&I_{H^2_{\clh}(\mathbb{D})}\end{bmatrix},\]that is, $N_{\clh}$
is an isometry. Moreover, one can check immediately that
\[N^k_{\clh} :=
\begin{bmatrix}M^k& 0\\ * &M^k_z\end{bmatrix}, \quad \quad \quad (k \in
\mathbb{N})\]which along with the isometric embedding $\Pi_N \in
\clb(\clh, \cln_{\clh})$ defined by $\Pi_{N} h = h \oplus 0$, for
all $h \in \clh$, implies that $\Pi_N$ is an isometric dilation of
$\clh$ with the isometric map $N_{\clh}$. Then $\tilde{\Pi}_{N} :=
P_{\tilde{\cln}_{\clh}}\Pi_{N} \in \clb(\clh, \tilde{\cln}_{\clh})$
is the minimal isometric dilation of $\clh$, where
$\tilde{\cln}_{\clh} = \overline{\mbox{span}} \{N_{\clh}^k \clh : k
\in \mathbb{N}\}$ and $\tilde{N}_{\clh} =
N_{\clh}|_{\tilde{\cln}_{\clh}}$. \qed

Sz.-Nagy's minimal isometric dilation is unique in the following
sense: if $\Pi \in \clb(\clh, \clm)$ is a minimal isometric dilation
of $\clh$ with isometry $V$, then there exists a (unique) unitary
$\Phi : \tilde{\cln}_{\clh} \raro \clm$ such that $V \Phi = \Phi
\tilde{N}_{\clh}$.

It is also worth mentioning that the Sch\"affer isometric dilation
of $\clh$ is always minimal. The \textit{Sch\"affer dilation space}
is defined by $\cls_{\clh} := \clh \oplus
H^2_{\cld_{\clh}}(\mathbb{D})$ with
\[S_{\clh} :=
\begin{bmatrix}M& 0\\\bm{D}_{\clh}&M_z\end{bmatrix},\]where $\cld_{\clh} =
\overline{\mbox{ran}}D_{\clh}$ (see \cite{NF}).

The von Neumann inequality \cite{vN} follows from the isometric, and
hence unitary (cf. \cite{NF}), dilation theorem for contractive
Hilbert modules.

\begin{Theorem}
Let $\clh$ be a Hilbert module over $\mathbb{C}[z]$. Then $\clh$ is
contractive if and only if \[\|p(M)\| \leq \|p\|_{\infty} =
\mathop{\mbox{max}}_{|z| \leq 1} |p(z)|. \quad \quad (p \in
\mathbb{C}[z])\]
\end{Theorem}
As a consequence, the polynomial functional calculus of a
contractive Hilbert module over $\mathbb{C}[z]$ extends to the disc
algebra $A(\mathbb{D})$, where \[A(\mathbb{D}) = \clo(\mathbb{D})
\cap \mathbb{C}(\overline{\mathbb{D}}) =
\overline{\mathbb{C}[z]}^{\|\cdot\|_{\infty}}.\]This implies that
\[\|f(M)\| \leq \|f\|_{\infty},\] for all $f \in A(\mathbb{D})$.

\subsection{Hilbert modules over $A(\Omega)$}\label{HMC}

Let $\Omega$ be a domain in $\mathbb{C}^n$ and $A(\Omega)$ be the
unital Banach algebra obtained from the closure in the supremum norm
on $\Omega$ of all functions holomorphic in some neighborhood of the
closure of $\Omega$. The most classical and familiar examples of
$A(\Omega)$ are the ball algebra $A(\mathbb{B}^n)$ and the polydisc
algebra $A(\mathbb{D}^n)$.

Now let $\clh$ be a Hilbert space and $\pi$ a norm continuous unital
algebra homomorphism from the Banach algebra $A(\Omega)$ to the
$C^*$-algebra $\clb(\clh)$. Then the Hilbert space $\clh$ is said to
be a \textit{Hilbert module over $A(\Omega)$} if it is an
$A(\Omega)$-module in the sense of algebra, \[A(\Omega) \times \clh
\raro \clh, \quad\quad (\varphi, f) \mapsto \varphi \cdot f =
\pi(\varphi)h,\]with the additional property that the module
multiplication $A(\Omega) \times \clh \raro \clh$ is norm
continuous. We say that a Hilbert module $\clh$ over $A(\Omega)$ is
contractive if $\pi$ is a contraction, that is, \[\|\varphi \cdot
f\|_{\clh} \leq \|\varphi\|_{A(\Omega)} \|f\|_{\clh}. \quad \quad
(\varphi \in A(\Omega), \, f \in \clh)\]

The following are some important and instructive examples of
contractive Hilbert modules.

(1) the \textit{Hardy module} $H^2(\mathbb{D}^n)$ (\cite{FS},
\cite{R}), the closure of $\mathbb{C}[\z]$ in $L^2(\mathbb{T}^n)$,
over $A(\mathbb{D}^n)$,

(2) the \textit{Hardy module} $H^2(\mathbb{B}^n)$ \cite{R-u}, the
closure of $\mathbb{C}[\z]$ in $L^2(\partial \mathbb{B}^n)$, over
$A(\mathbb{B}^n)$ and

(3) the \textit{Bergman module over the ball} $L^2_a(\mathbb{B}^n)$,
the closure of $A(\mathbb{B}^n)$ in $L^2(\mathbb{B}^n)$, over
$A(\mathbb{B}^n)$.

(4) Quotient modules and submodules of (1), (2) and (3) over the
corresponding algebras.

\subsection{Module tensor products and localizations}
Module tensor product and localizations are in the center of
commutative algebra and algebraic geometry. The notion of module
tensor products and localizations for Hilbert modules introduced by
Douglas and Paulsen \cite{DP} was one of the inspiration points for
Hilbert module method to operator theory.

Let $\clh_1$ and $\clh_2$ be two Hilbert modules over $A(\Omega)$
and $\clh_1 \otimes \clh_2$ be the Hilbert space tensor product.
Then $\clh_1 \otimes \clh_2$ turns into both a left and a right
$A(\Omega)$ modules, $A(\Omega) \times \clh_1 \otimes \clh_2 \raro
\clh_1 \otimes \clh_2$, by setting
\[(\varphi, h_1 \otimes h_2) \mapsto (\varphi \cdot h_1)
\otimes h_2, \quad \mbox{and} \quad (\varphi, h_1 \otimes h_2)
\mapsto h_1 \otimes (\varphi \cdot h_2),\]respectively. Note
that\[\cln = \overline{\mbox{span}}\{(\varphi \cdot h_1) \otimes h_2
- h_1 \otimes (\varphi \cdot h_2): h_1 \in \clh_1, h_2 \in \clh_2,
\varphi \in A(\Omega)\}\]is both a left and a right
$A(\Omega)$-submodule of $\clh_1 \otimes \clh_2$. Then $\cln^{\perp}
(\cong (\clh_1 \otimes \clh_2)/\cln)$ is both a left and a right
$A(\Omega)$-quotient module and
\[P_{\cln^{\perp}} ((\varphi \cdot h_1) \otimes h_2) =
P_{\cln^{\perp}} (h_1 \otimes (\varphi \cdot h_2)),\]for all $h_1
\in \clh_1, h_2 \in \clh_2$ and $\varphi \in A(\Omega)$. In
conclusion, these quotient modules are isomorphic Hilbert modules
over $A(\Omega)$, which will denote by $\clh_1 \otimes_{A(\Omega)}
\clh_2$ and referred as the \textit{module tensor product} of the
Hilbert modules $\clh_1$ and $\clh_2$ over $A(\Omega)$.

For each $\w \in \Omega$ denote by $\mathbb{C}_{\w}$ the one
dimensional Hilbert module over $A(\Omega)$: \[A(\Omega) \times
\mathbb{C}_{\w} \raro \mathbb{C}_{\w}, \quad \quad (\varphi,
\lambda) \mapsto \varphi(\w) \lambda.\] Further, for each $\w \in
\Omega$ denote by $A(\Omega)_{\w}$ the set of functions in
$A(\Omega)$ vanishing at $\w$, that is,
\[A(\Omega)_{\w} = \{\varphi \in A(\Omega) : \varphi(\w) = 0\}.\]

Let $\clh$ be a Hilbert module over $A(\Omega)$ and $\w \in \Omega$.
Then the module tensor product $\clh \otimes_{A(\Omega)}
\mathbb{C}_{\w}$ is called the \textit{localization} of the Hilbert
module $\clh$ at $\w$.

It is easy to see that \[\clh_{\w} : = \overline{\mbox{span}}\{
\varphi f : \varphi \in A(\Omega)_{\w}, f \in \clh\}.\]is a
submodule of $\clh$ for each $\w \in \Omega$. Moreover, the quotient
module $\clh/ \clh_{\w}$ is canonically isomorphic to $\clh
\otimes_{A(\Omega)} \mathbb{C}_{\w}$, the localization of $\clh$ at
$\w \in \Omega$, in the following sense
\[\clh \otimes_{A(\Omega)} \mathbb{C}_{\w} \raro \clh/ \clh_{\w},
\quad \quad f \otimes_{A(\Omega)} 1 \mapsto P_{\clh/ \clh_{\w}} f.\]

The following list of examples of localizations will be useful in a
number of occasions later.

\NI \textsf{(1)} For all $\w \in \mathbb{D}^n$, \[H^2(\mathbb{D}^n)
\otimes_{A(\mathbb{D}^n)} \mathbb{C}_{\w} \cong
\mathbb{C}_{\w},\]where $\cong$ stands for module isomorphism.

\NI\textsf{(2)} Let $\clh = H^2(\mathbb{B}^n)$ or
$L^2_a(\mathbb{B}^n)$. Then for all $\w \in \mathbb{B}^n$,
\[\clh \otimes_{A(\mathbb{B}^n)} \mathbb{C}_{\w} \cong
\mathbb{C}_{\w}.\]

\NI\textsf{(3)} Let $H^2(\mathbb{D}^2)_{{0}} = \{f \in
H^2(\mathbb{D}^2): f({0}) = 0\}$, the submodule of
$H^2(\mathbb{D}^2)$ of functions vanishing at the origin. Then \[
H^2(\mathbb{D}^2)_{{0}} \otimes_{A(\mathbb{D}^2)} \mathbb{C}_{\w} =
\left\{
\begin{array}{ll}
\mathbb{C}_{\w} & \;\;\mbox{if $\w \neq 0$};\\
\mathbb{C}_0 \oplus \mathbb{C}_0 & \;\;\mbox{if $\w =
0$}.\end{array} \right. \]

\NI \textbf{Further results and comments:}

\begin{enumerate}

\item von Neumann inequality says that one can extend the functional
calculus from $\mathbb{C}[z]$ to $A(\mathbb{D})$ for contractive
Hilbert module over $\mathbb{C}[z]$. Another approach to extend the
functional calculus is to consider the rational functions. More
precisely, let $K$ be a non-empty compact subset of $\mathbb{C}$ and
$T \in \clb(\clh)$. Denote $Rat(X)$ the set of rational functions
with poles off $K$. Then $K$ is a \textit{spectral set} for $T$ if
$\sigma(T) \subseteq K$ and
\[\|f(T)\| \leq \|f\|_K := \mbox{sup~}\{|f(z)| : z \in K\}.\]
The notion of spectral set was introduced by J. von Neumann in
\cite{vN} where he proved that the closed unit disk is a spectral
set of a bounded linear operator on a Hilbert space if and only if
the operator is a contraction. Also recall that a bounded linear
operator $T$ on $\clh$ has a normal $\partial K$-dilation if there
exists a normal operator $N$ on $\clk \supseteq \clh$ such that
$\sigma(N) \subseteq
\partial K$ and
\[P_{\clh} f(N)|_{\clh} = f(T). \quad \quad (f \in \,Rat(K))\]The Sz. Nagy
dilation theory shows that every contraction has a normal $\partial
\overline{\mathbb{D}}$-dilation. It is known that the normal
$\partial K$-dilation holds if $K$ is the closure of an annulus
\cite{A-Annals} and fails, in general, when $K$ is a triply
connected domain in $\mathbb{C}$ \cite{DM} (see also \cite{AHR},
\cite{AY} and \cite{S-SB}).

\item Ando's theorem \cite{An} extends Sz.-Nagy's unitary dilation
result to a pair of operators, that is, any pair of commuting
contractions has a unitary dilation. However, the Ando dilation is
not unique and it fails for three or more operators (see \cite{Pa}
and \cite{Va}).

\item The Ando dilation theorem is closely related to the
commutant lifting theorem (see \cite{Pa}).

\item A $2$-variables analogue of von Neumann's inequality follows from Ando's dilation theorem \cite{An}.
It is well known that for $n$-tuples of operators, $n \geq 3$, the
von Neumann inequality fails in general. In \cite{AKVW}, Anatolii,
Kaliuzhnyi-Verbovetskyi, Vinnikov and Woerdeman proved a several
variables analogue of von Neumann's inequality for a class of
commuting $n$-tuples of strict contractions.

\item The notion of module tensor product is due to Douglas (see
\cite{DP}).\

\item The notion of localization of Hilbert modules, however, is by far not enough. The
computation of higher order localizations is another important issue
in the theory of Hilbert modules over $\mathbb{C}[\z]$, which in
general can be very difficult \cite{ChenD}.

\item In connection with this section see \cite{DP}, \cite{NF},
\cite{S-SB}.
\end{enumerate}

\section{Hilbert modules of holomorphic functions}\label{HMHF}

In various parts of operator theory and functional analysis, one is
confronted with Hilbert spaces of functions, such that it is both
simple and instructive to deal with a large class of operators (cf.
\cite{AE}). The purpose of this section is to provide a brief
introduction to the theory of Hilbert modules of holomorphic
functions that will be used in subsequent sections.

\subsection{Reproducing kernel Hilbert modules}\label{sub-RKHM}
A natural source of Hilbert module comes from the study of
reproducing kernel Hilbert spaces (cf. \cite{A}, \cite{AM-book},
\cite{CS}) on domains in $\mathbb{C}^n$.

Let $X$ be a non-empty set, and $\cle$ a Hilbert space. An
operator-valued function $K: X \times X \rightarrow \clb(\cle)$ is
said to be \textit{positive definite kernel} if \[ \sum_{i, j =
1}^{k} \langle K(z_i, z_j) \eta_j, \eta_i \rangle \geq 0,\] for all
$\eta_i \in \cle, \, z_i \in X$, $i=1, \cdots, k$, and $k \in
\mathbb{N}$. Given such a positive definite kernel $K$ on $X$, let
$\clh_K$ be the Hilbert space completion of the linear span of the
set $\{K(\cdot, w) \eta : w \in X, \eta \in \cle\}$ with respect to
the inner product
\[\langle K(\cdot, w) \eta, K(\cdot, z) \zeta\rangle_{\clh_K} = \langle K(z,
w) \eta, \zeta\rangle_{\cle},\]for all $z, w \in X$ and $\eta, \zeta
\in \cle$. Therefore, $\clh_K$ is a Hilbert space of $\cle$-valued
functions on $X$. The kernel function $K$ has the reproducing
property:
\[\langle f, K(\cdot, z) \eta\rangle_{\clh_K} = \langle f(z),
\eta\rangle_{\cle},\]for all $z \in X$, $f \in \clh_K$ and $\eta \in
\cle$. In particular, for each $z \in X$, the evaluation operator
${ev}_z : \clh_K \rightarrow \cle$ defined by \[\langle {ev}_z(f),
\eta\rangle_{\cle} = \langle f, K(\cdot, z) \eta\rangle_{\clh_K},
\quad \quad (\eta \in \cle, f \in \clh_K)\] is bounded. Conversely,
let $\clh$ be a Hilbert space of functions from $X$ to $\cle$ with
bounded and non-zero evaluation operators $ev_z$ for all $z \in X$.
Therefore, $\clh$ is a reproducing kernel Hilbert space with
reproducing kernel
\[K(z, w) = {ev}_z \circ {ev}_w^* \in \clb(\cle). \quad \quad (z, w
\in X)\]

\NI Now let $X = \Omega$ a domain in $\mathbb{C}^n$ and $K : \Omega
\times \Omega \rightarrow \clb(\cle)$ be a kernel function,
holomorphic in the first variable and anti-holomorphic in the second
variable. Then $\clh_K$ is a Hilbert space of holomorphic functions
on $\Omega$ (cf. \cite{CS}).

A Hilbert module $\clh_K$ is said to be \textit{reproducing kernel
Hilbert module} over $\Omega$ if $\clh_K \subseteq \clo(\Omega,
\cle)$ and for each $1 \leq i \leq n$,
\[M_i f = z_i f,\]where \[(z_i f)(\w) = w_i f(\w). \quad \quad \quad \quad \quad (f \in \clh_K, \w \in
\Omega)\]It is easy to verify that
\[M_{z_i}^* (K(\cdot, \w) \eta) = \bar{w}_i K(\cdot, \w) \eta,\]for
all $\w \in \Omega, \eta \in \cle$ and $i=1, \ldots, n$.

\textit{In most of the following the module maps $\{M_i\}_{i=1}^n$
of a reproducing kernel Hilbert module will simply be denoted by the
multiplication operators $\{M_{z_i}\}_{i=1}^n$ by the coordinate
functions $\{{z_i}\}_{i=1}^n$.}

\NI \textsf{Examples:} (1) The Drury-Arveson module, denoted by
$H^2_n$, is the reproducing kernel Hilbert module corresponding to
the kernel $k_n : \mathbb{B}^n \times \mathbb{B}^n \raro
\mathbb{C}$, where \[k_n(\z, \w) = (1 - \sum_{i=1}^n z_i
\bar{w}_i)^{-1}. \quad \quad (\z, \w \in \mathbb{B}^n)\]

\NI (2) Suppose $\alpha > n$. The weighted Bergman space $L^2_{a,
\alpha}(\mathbb{B}^n)$ (see \cite{Zhu-Z}) is a reproducing kernel
Hilbert space with kernel function \[k_{\alpha}(\z, \w) =
\frac{1}{(1 - \langle \z, \w \rangle_{\mathbb{C}^n})^{\alpha}}.
\quad \quad (\z, \w \in \mathbb{B}^n)\]When $\alpha = n$, $L^2_{a,
\alpha}(\mathbb{B}^n)$ is the usual Hardy module
$H^2(\mathbb{B}^n)$.

\NI (3) The kernel function for the Dirichlet module (see
\cite{Zhu-TAMS}) $\cld(\mathbb{B}^n)$ is given by
\[k_{\cld(\mathbb{B}^n)}(\z, \w) = 1 + \log\frac{1}{1 - \langle \z,
\w \rangle_{\mathbb{C}^n}}. \quad \quad (\z, \w \in \mathbb{B}^n)\]

\NI (4) $H^2(\mathbb{D}^n)$, the Hardy module over $\mathbb{D}^n$,
is given by the reproducing kernel \[\mathbb{S}_n(\z, \w) =
\mathop{\Pi}_{i=1}^n (1 - z_i \bar{w}_i)^{-1}. \quad \quad (\z, \w
\in \mathbb{D}^n)\]

Finally, let $I$ be a non-empty set and \[l^2(I) = \{f : I \raro
\mathbb{C} : \mathop{\sum}_{i \in I} |f(i)|^2 < \infty\}.\]Then
$l^2(I)$ is a reproducing kernel Hilbert space with kernel $k(i, j)
= \delta_{ij}$ for all $(i, j) \in I \times I$. Moreover,
$\{k(\cdot, j) : j \in I\}$ is an orthonormal basis of $l^2(I)$. In
general, $l^2(I)$ is not a reproducing kernel Hilbert module.

\subsection{Cowen-Douglas Hilbert modules}

Let $m$ be a positive integer. A class of Hilbert modules over
$\Omega \subseteq \mathbb{C}$, denoted by $B_m(\Omega)$, was
introduced by Cowen and Douglas in \cite{CD}. This notion was
extended to the multivariable setting, for $\Omega \subseteq
\mathbb{C}^n$, by Curto and Salinas \cite{CS} and by Chen and
Douglas \cite{ChenD}. See also \cite{CD2}.

\begin{Definition}\label{CD-Defn}
Let $\Omega$ be a domain in $\mathbb{C}^n$ and $m$ be a positive
integer. Then a Hilbert module $\clh$ over $\mathbb{C}[\z]$ is said
to be in $B_m^*(\Omega)$ if

(i) the column operator $(M - \w I_{\clh})^* : \clh \rightarrow
\clh^n$ defined by

\[(M - \w I_{\clh})^* h = (M_{1} - w_1 I_{\clh})^* h \oplus
\cdots \oplus (M_{n} - w_n I_{\clh})^* h,\quad \quad (h \in \clh)\]
has closed range for all $\w \in \Omega$, where $\clh^n = \clh
\oplus \cdots \oplus \clh$.

(ii) $\mbox{dim~} \mbox{ker~}(M - \w I_{\clh})^* = \mbox{dim}
[~\cap_{i=1}^n \mbox{ker} (M_{i} - w_i I_{\clh})^*] = m$ for all $\w
\in \Omega$, and

(iii) $\bigvee_{{\w} \in \Omega} \mbox{ker~}(M - \w I_{\clh})^* =
\clh$.
\end{Definition}

Given a Hilbert module $\clh$ in $B_m^*(\Omega)$, define
\[E^*_{\clh} = \mathop{\bigcup}_{{\w} \in \Omega} \{\bar{\w}\}
\times \mbox{ker} (M - \w I_{\clh})^*.\] Then the mapping $\w
\mapsto E^*_{\clh}(\w) := \{\bar{\w}\} \times \mbox{ker~} (M - \w
I_{\clh})$ defines a rank $m$ hermitian anti-holomorphic vector
bundle over $\Omega$. For a proof of this fact, the reader is
referred to \cite{CD}, \cite{CD2}, \cite{CS} and \cite{ES}.

The fundamental relation between $\clh \in B^*_m(\Omega)$ and the
associated anti-holomorphic hermitian vector bundle \cite{W} over
$\Omega$ defined by

\[\begin{CD}E^*_{\clh} : \mbox{ker~}(M - \w I_{\clh})^* \\ \downarrow\\
\w
\end{CD}\]is the following identification:

\begin{Theorem}\label{CD-rigid}
Let $\Omega = \mathbb{B}^n$ or $\mathbb{D}^n$ and $\clh,
\tilde{\clh} \in B^*_m(\Omega)$. Then $\clh \cong \tilde{\clh}$ if
and only if the complex bundles $E_{\clh}^*$ and
$E^*_{\tilde{\clh}}$ are equivalent as Hermitian anti-holomorphic
vector bundles.
\end{Theorem}

Note that for $U$ an open subset of $\Omega$, the anti-holomorphic
sections of $E^*_{\clh}$ over $U$ are given by $\gamma_f : U \raro
E^*_{\clh}$, where $\gamma_f(\w) = (\bar{\w}, f(\w))$ and $f : U
\raro \clh$ is an anti-holomorphic function with $f(\w) \in \ker (M
- \w I_{\clh})^*$ for all $\w \in U$.

The Grauert's theorem asserts that the anti-holomorphic vector
bundle $E^*_{\clh}$ over a domain in $\mathbb{C}$ or a contractible
domain of holomorphy in $\mathbb{C}^n$ is holomorphically trivial,
that is, $E^*_{\clh}$ possesses a global anti-holomorphic frame. In
particular, there exists anti-holomorphic functions $\{s_i\}_{i=1}^m
\subseteq \clo^*(\Omega, \clh)$ such that $\{s_i(\w)\}_{i=1}^m$ is a
basis of $\mbox{ker~} (M - \w I_{\clh})$ for all $\w \in \Omega$.
Moreover, $\clh$ is unitarily equivalent to a reproducing kernel
Hilbert module with $\clb(\mathbb{C}^m)$-valued kernel (see
\cite{Al}, \cite{CS}, \cite{ES}).

\begin{Theorem}
Let $\clh \in B^*_m(\Omega)$ where $\Omega$ be a domain in
$\mathbb{C}$ or a contractible domain of holomorphy in
$\mathbb{C}^n$. Then there exits a reproducing kernel Hilbert module
$\clh_K \subseteq \clo(\Omega, \mathbb{C}^m)$ such that $\clh \cong
\clh_K$.
\end{Theorem}

\NI \textsf{Proof.} Define $J_s : \clh \raro \clo(\Omega,
\mathbb{C}^m)$ by
\[(J_s(f))(\w) = (\langle f, s_1(\w) \rangle_{\clh}, \ldots, \langle f,
s_m(\w)\rangle_{\clh}). \quad \quad (f \in \clh, \w \in
\Omega)\]Note that $J_s$ is an injective map. Consequently, the
space $\clh_{J_s} := \mbox{ran} {J_s} \subseteq \clo(\Omega,
\mathbb{C}^m)$ equipped with the norm \[\|J_s f\|_{\clh_{J_s}} :=
\|f\|_{\clh}, \quad \quad (f \in \clh)\]is a $\mathbb{C}^m$-valued
reproducing kernel Hilbert space with kernel $K_s : \Omega \times
\Omega \raro \clb(\mathbb{C}^m)$ given by the "Gram matrix" of the
frame $\{s_i(\w) : 1 \leq i \leq m \}$:
\[K_s(\z, \w) = \big( \langle s_j(\w), s_i(\z)\rangle_{\clh})_{i,j=1}^m.
\quad \quad (\z, \w \in \Omega)\] Further, note that
\[\begin{split}(J_s M_i f)(\w) & =  (\langle M_i f, s_1(\w) \rangle_{\clh}),
\ldots, \langle M_i f, s_m(\w)\rangle_{\clh})) \\& = (\langle f,
M_i^* s_1(\w) \rangle_{\clh}), \ldots, \langle f, M_i^*
s_m(\w)\rangle_{\clh}))\\& = w_i(\langle f, s_1(\w) \rangle_{\clh}),
\ldots, \langle f, s_m(\w)\rangle_{\clh}))
\\& = (M_{z_i} J_s f)(\w),
\end{split}\] for all $f \in
\clh$ and $\w \in \Omega$. This implies that $J_s M_i = M_{z_i} J_s$
for all $1 \leq i \leq n$ and hence the Hilbert module $\clh$ is
module isometric isomorphic with the reproducing kernel Hilbert
module $\clh_{J_s}$. \qed

If $E^*_{\clh}$ is not trivial, then we can use an anti-holomorphic
frame over an open subset $U \subseteq \Omega$ to define a kernel
function $K_U$ on $U$. Since a domain is connected, one can show
that $\clh_{K_U} \cong \clh$. One way to obtain a local frame is to
identify the fiber of the dual vector bundle $E_{\clh}$ with
$\clh/I_{\w} \cdot \clh \cong \mathbb{C}^m \cong \mbox{span}
\{s_i(\w) : 1 \leq i \leq m\}$, where $I_{\w} = \{ p \in
\mathbb{C}[\z] : p(\w) = 0\}$ is the maximal ideal of
$\mathbb{C}[\z]$ at $\w \in \Omega$.

The curvature of the bundle $E^*_{\clh}$ for the Chern connection
determined by the metric defined by the Gram matrix or, if
$E^*_{\clh}$ is not trivial, then with the inner product on
$E^*_{\clh}(\w) = \mbox{ker}(M_z - \w I_{\clh})^* \subseteq \clh$,
is given by
\[\clk_{E^*_{\clh}} (\w) = (\bar{\partial}_j \{K(\w, \w)^{-1} \partial_i K(\w, \w)\})_{i,j=1}^n,\]
for all $\w \in \Omega$. Note that the representation of the
curvature matrix defined above is with respect to the basis of
two-forms $\{dw_i \wedge d\bar{w}_j : 1 \leq i,j \leq n\}$. In
particular, for a line bundle, that is, when $m=1$, the curvature
form is given by
\[\begin{split}\clk_{E^*_{\clh}}(\w) & = \bar{\partial} K(\w,\w)^{-1}
\partial K(\w, \w) = -  \partial \bar{\partial} \log \|K(\cdot, \w)\|^2
\\ & = - \sum_{i, j =
1}^n \frac{\partial^2}{\partial w_i \partial\bar{w_j}} \log K(\w,\w)
dw_i \wedge d\bar{w}_j. \quad \quad (\w \in \Omega)\end{split}\]

The Hardy modules $H^2(\mathbb{B}^n)$ and $H^2(\mathbb{D}^n)$, the
Bergman modules $L^2_a(\mathbb{B}^n)$ and $L^2_a(\mathbb{D}^n)$, the
weighted Bergman modules $L^2_{a, \alpha}(\mathbb{B}^n)$ ($\alpha
>n$) and the Drury-Arveson module $H^2_n$ are the standard examples
of Hilbert modules in $B^*_1(\Omega)$ with $\Omega = \mathbb{B}^n$
or $\mathbb{D}^n$. A further source of Hilbert modules in
$B^*_m(\Omega)$ is a family of some quotient Hilbert modules, where
the standard examples are used as building blocks (see Section 2 in
\cite{S-HM}).

\subsection{Quasi-free Hilbert modules}

Besides reproducing kernel Hilbert modules, there is another class
of function Hilbert spaces which will be frequently used throughout
this article. These are the quasi-free Hilbert modules.

Recall that the Hardy module and the weighted Bergman modules over
$\mathbb{D}^n$ (or $\mathbb{B}^n$) are singly-generated Hilbert
module over $A(\mathbb{D}^n)$ (over $A(\mathbb{B}^n)$). In other
words, these modules are the Hilbert space completion of
$A(\Omega)$. More generally, every cyclic or singly-generated
bounded Hilbert module over $\mathcal A(\Omega)$ is obtained as a
Hilbert space completion of $\mathcal A(\Omega)$.

On the other hand, finitely generated free modules over $A(\Omega)$,
in the sense of commutative algebra, have the form $A(\Omega)
\otimes_{alg} l^2_m$ for some $m \in \mathbb{N}$ (see \cite{DE}).
However, the algebraic tensor product $A(\Omega) \otimes_{alg}
l^2_m$ is not a Hilbert space. In order to construct ``free Hilbert
modules'' we consider Hilbert space completions of free modules
$A(\Omega) \otimes_{alg} l^2_m$:

\NI Let $m \geq 1$. A Hilbert space $\clr$ is said to be
\textit{quasi-free Hilbert module over $A(\Omega)$ and of rank $m$}
if $\clr$ is a Hilbert space completion of the algebraic tensor
product $A(\Omega) \otimes_{alg} l^2_m$ and
\begin{enumerate}
\item  multiplication by functions in $\mathcal A(\Omega)$ define bounded operators on $\clr$,
\item the evaluation operators $ev_{\w} : \clr \to l_m^2$ are locally uniformly
bounded on $\Omega$, and
\item  a sequence $\{f_k\} \subseteq \mathcal A(\Omega) \otimes l_m^2$ that is
Cauchy in the norm of $\clr$ converges to $0$ in the norm of $\clr$
if and only if $ev_{\w}(f_k)$ converges to $0$ in $l^2_m$ for $\w
\in \Omega$.
\end{enumerate}

Condition (1) implies that $\clr$ is a bounded Hilbert module over
$A(\Omega)$. Condition (2) ensures that $\clr$ can be identified
with a Hilbert space of $l^2_m$-valued holomorphic functions on
$\Omega$ and condition (3) implies that the limit function of a
Cauchy sequence in $A(\Omega) \otimes_{alg} l^2_m$ vanishes
identically if and only if the limit in the $\clr$-norm is the zero
function. In other words, a quasi-free Hilbert module $\clr$ over
$A(\Omega)$ is a finitely generated reproducing kernel Hilbert
module where the kernel function $K: \Omega\times \Omega \raro
\clb(l_m^2)$ is holomorphic in the first variable and
anti-holomorphic in the second variable.

In some instances, such as the Drury-Arveson module $H^2_n$, this
definition does not apply. In such cases $\clr$ is defined to be the
completion of the polynomial algebra $\mathbb{C}[\z]$ relative to an
inner product on it assuming that each $p(\z)$ in
$\mathbb{C}[\bm{z}]$ defines a bounded operator on $\clr$ but there
is no uniform bound. Hence, in this case $\clr$ is a Hilbert module
over $\mathbb{C}[\bm{z}]$.

\subsection{Multipliers}

Given $\cle$- and $\cle_*$-valued reproducing kernel Hilbert modules
$\clh$ and $\clh_*$, respectively, over $\Omega$, a function
$\varphi : \Omega \rightarrow \clb(\cle, \cle_*)$ is said to be a
\textit{multiplier} if $\varphi f \in \clh_*$, where $(\varphi
f)(\w) = \varphi(\w) f(\w)$ for $f \in \clh$ and $\w \in \Omega$.
The set of all such multipliers is denoted by $\clm(\clh, \clh_*)$
or simply $\clm$ if $\clh$ and $\clh_*$ are clear from the context
(cf. \cite{VB}). By the closed graph theorem, each $\varphi \in
\clm(\clh, \clh_*)$ induces a bounded linear map $M_{\varphi} : \clh
\rightarrow \clh_*$ (cf. \cite{H-book}) defined by
\[M_{\varphi} f = \varphi f,\]for all $f \in \clh_K$. Consequently,
$\clm(\clh, \clh_*)$ is a Banach space with
\[\|\varphi\|_{\clm(\clh, \clh_*)} = \|M_{\varphi}\|_{\clb(\clh,\clh_*)}.\]
\NI For $\clh = \clh_*$, $\clm(\clh) = \clm(\clh, \clh)$ is a Banach
algebra with this norm.

Let $\clr \subseteq \clo(\Omega, \mathbb{C})$ be a reproducing
kernel Hilbert module with kernel $k_{\clr}$ and $\cle$ be a Hilbert
space. Then $\clr \otimes \cle$ is a reproducing kernel Hilbert
module with kernel function $(\z, \w) \mapsto k_{\clr}(\z, \w)
I_{\cle}$. By $\clm_{\clb(\cle, \cle_*)}(\clr)$ we denote the set of
all multipliers $\clm(\clr \otimes \cle, \clr \otimes \cle_*)$.

The following characterization result is well known and easy to
prove.

\begin{Theorem}
Let $X$ be a non-empty set and for $i = 1, 2$, $K_i : X \times X
\raro \clb(\cle_i)$ be positive definite kernel functions with
reproducing kernel Hilbert spaces $\clh_{K_i}$. Suppose also that
$\Theta : X \raro \clb(\cle_1, \cle_2)$ is a function. Then the
following are equivalent:

(1) $\Theta \in \clm (\clh_{K_1}, \clh_{K_2})$.

(2) There exists a constant $c > 0$ such that \[(x, y) \raro c^2
K_2(x, y) - \Theta(x) K_1(x, y) \Theta(y)^*\]is positive definite.
In this case, the multiplier norm of $\Theta$ is the infimum of all
such constants $c >0$. Moreover, the infimum is achieved.
\end{Theorem}

\NI\textsf{Examples:}

(1) For the Drury-Aveson space $H^2_n$, the multiplier space is
given by \[\clm_{\clb(\cle, \cle_*)}(H^2_n) = \{\Theta \in
\clo(\mathbb{B}^n, \clb(\cle, \cle_*)) : \mbox{sup} \|\Theta(rT)\| <
\infty\},\]where the supremum ranges over $0 < r < 1$ and commuting
$n$-tuples $(T_, \ldots, T_n)$ on Hilbert spaces $\clh$ such that
$\sum_{i=1}^n T_i T^*_i \leq I_{\clh}$ (see \cite{EP02}, \cite{BTV}
for more details).

(2) Let $\clh = H^2(\mathbb{B}^n)$ or $L^2_a(\mathbb{B}^n)$. Then
\[\clm_{\clb(\cle, \cle_*)}(\clh) = H^\infty_{\clb(\cle,
\cle_*)}(\mathbb{B}^n).\]

(3) Let $\clh = H^2(\mathbb{D}^n)$ or $L^2_a(\mathbb{D}^n)$. Then
\[\clm_{\clb(\cle, \cle_*)}(\clh) = H^\infty_{\clb(\cle,
\cle_*)}(\mathbb{D}^n).\]

One striking fact about the Dirichlet space is that the multiplier
space $\clm(\cld(\mathbb{D}))$ is a proper subset of
$H^\infty(\mathbb{D})$ (see \cite{St}). Also, it is bounded but not
a contractive Hilbert module over $\mathbb{C}[z]$. Note that also
the multiplier space $\clm(H^2_n)$ is a proper subspace of
$H^\infty(\mathbb{B}^n)$. Moreover, $\clm(H^2_n)$ does not contain
the ball algebra $A(\mathbb{B}^n)$ (see \cite{D-78}, \cite{AAM}).

This subsection concludes with a definition. Let $\Theta_i \in
\clm_{\clb(\cle_i, \cle_{*i})}(\clr)$ and $i = 1, 2$. Then
$\Theta_1$ and $\Theta_2$ are said to \textit{coincide}, denoted by
$\Theta_1 \cong \Theta_2$, if there exists unitary operators $\tau :
\cle_1\raro \cle_{2}$ and $\tau_* : \cle_{*1} \raro \cle_{*2}$ such
that the following diagram commutes:
\[\begin{CD}\clr \otimes \cle_1 @>M_{\Theta_1} >> \clr \otimes
\cle_{*1}\\ @ V I_{\clr} \otimes \tau VV @VI_{\clr} \otimes \tau_* VV\\\
\clr \otimes \cle_2 @>M_{\Theta_2}>> \clr \otimes \cle_{*2}
\end{CD}\]

\NI \textbf{Further results and comments:}

\begin{enumerate}

\item Let $\Omega \subseteq \mathbb{C}$ and $V^*(\clh)$ be the von Neumann
algebra of operators commuting with both $M_z$ and $M_z^*$. Note
that projections in $V^*(\clh)$, or reducing submodules of $\clh$,
are in one-to-one correspondence with reducing subbundles of
$E^*_{\clh}$. A subbundle $F$ of an anti-holomorphic Hermitian
vector bundle $E$ is said to be a \textit{reducing subbundle} if
both $F$ and its orthogonal complement $F^{\perp}$ in E are
anti-holomorphic subbundles.

\NI Also note that if $S$ is an operator commuting with $M_z^*$,
then $S E_{\clh}^*(w) \subseteq E^*_{\clh}(w)$ for each $w \in
\Omega$ and hence $S$ induces a holomorphic bundle map, denoted by
$\Gamma(S)$, on $E^*_{\clh}$. In \cite{CDG}, Chen, Douglas and Guo
proved that if $S$ lies in $V^*(\clh)$, then $\Gamma(S)$ is not only
anti-holomorphic, but also connection-preserving.

\begin{Theorem} Let $\clh \in B^*_m(\Omega)$ and $\Phi$
be a bundle map on $E^*_{\clh}$. There exists an operator $T_{\Phi}$
in $V^*(\clh)$ such that $\Phi = \Gamma(T_{\Phi})$ if and only if
$\Phi$ is connection preserving. Consequently, the map $\Gamma$ is a
$*$-isomorphism from $V^*(\clh)$ to connection-preserving bundle
maps on $E^*_{\clh}$.
\end{Theorem}

\item Let $\clh_1 \in B_{m_1}^*(\Omega)$ and $\clh_2 \in
B_{m_2}^*(\Omega)$ and $\Omega \subseteq \mathbb{C}$. It is natural
to ask the following question: Determine the Hilbert module $\clh$,
if such exists, in $B_{m_1 m_2}(\Omega)$ corresponding to the
anti-holomorphic vector bundle $E^*_{\clh_1} \otimes E^*_{\clh_2}$.
That is, find $\clh \in B_{m_1 m_2}(\Omega)$ such that $E^*_{\clh}
\cong E^*_{\clh_1} \otimes E^*_{\clh_2}$, where the equivalence is
in terms of the anti-holomorphic vector bundle isomorphism. In
\cite{QL}, Q. Lin proved the following remarkable result.

\begin{Theorem}
Let $\clh_1 \in B_{m_1}^*(\Omega)$ and $\clh_2 \in
B_{m_2}^*(\Omega)$ and $\Omega \subseteq \mathbb{C}$. Define \[\clh
= \mathop{\bigvee}_{z \in \Omega} [\mbox{ker}(M - z I_{\clh_1})^*
\otimes \mbox{ker}(M - z I_{\clh_2})^*]. \]Then $\clh$ is a
submodule of $\clh_1 \otimes \clh_2$, and the module multiplications
on $\clh$ coincides: $(M \otimes I_{\clh_2})|_{\clh} = (I_{\clh_1}
\otimes M)|_{\clh}$. Moreover, $\clh \in B^*_{m_1 m_2}(\Omega)$ and
$E^*_{\clh} \cong E^*_{\clh_1} \otimes E^*_{\clh_2}$.
\end{Theorem}

\item
In \cite{Z-00}, Zhu suggested an alternative approach to the
Cowen-Douglas theory based on the notion of spanning holomorphic
cross-sections. More precisely, let $\Omega \subseteq \mathbb{C}$
and $\clh \in B_m^*(\Omega)$. Then $E^*_{\clh}$ possesses a spanning
anti-holomorphic cross-section, that is, there is an
anti-holomorphic function $\gamma : \Omega \raro \clh$ such that
$\gamma(w) \in \mbox{ker~}(M - w I_{\clh})^*$ for all $w \in \Omega$
and $\clh$ is the closed linear span of the range of $\gamma$.

\NI More recently, Eschmeier and Schmitt \cite{ES} extended Zhu's
results to general domains in $\mathbb{C}^n$.

\item The concept of quasi-free Hilbert module is due to Douglas and Misra \cite{DM1},
\cite{DM2}. The notion is closely related to the generalized Bergman
kernel introduced by Curto and Salinas \cite{CS}.

\item For a systematic exposition of the theory of quasi-free Hilbert
modules, see the work by Chen \cite{LChen}.

\item In connection with Cowen-Douglas theory see Apostol and Martin
\cite{AM-CD}, McCarthy \cite{Mc} and Martin \cite{Mircea}.

\item In \cite{BJOT}, Barbian proved that an operator $T$ between
reproducing kernel Hilbert spaces is a multiplier if and only if $(T
f)(x) = 0$ holds for all $f$ and $x$ satisfying $f(x) = 0$.

\item The reader is referred to \cite{A}, \cite{AM-book}, \cite{CS}, \cite{DMV} and \cite{BM-book}
for some introduction to the general theory of reproducing kernel
Hilbert spaces. For recent results on reproducing kernel Hilbert
spaces see \cite{BJOT}, \cite{BIEOT} and the reference therein.

\end{enumerate}

\section{Contractive Hilbert modules over $A(\mathbb{D})$}\label{CHM}

This section gives a brief review of contractive Hilbert modules
over $A(\mathbb{D})$ and begins with the definition of free
resolutions from commutative algebra. The following subsection
recast the canonical model of Sz.-Nagy and Foias in terms of Hilbert
modules. It is proved that for a contractive Hilbert module over
$A(\mathbb{D})$ there exists a unique free resolution. The final
subsection is devoted to prove that the free resolutions of
contractive Hilbert modules over $A(\mathbb{D})$ are uniquely
determined by a nice class of bounded holomorphic functions on
$\mathbb{D}$.

\subsection{Free resolutions} The purpose of this subsection is to recall the notion of
\textit{free modules} from commutative algebra. Let $M$ be a module
over a commutative ring $R$. Then $M$ is free if and only if $M$ is
a direct sum of isomorphic copies of the underlying ring $R$.

It is well known and easy to see that every module has a free
resolution with respect to the underlying ring. More precisely,
given a module $M$ over a ring $R$, there exists a sequence of free
$R$-modules $\{F_i\}_{i=0}^{\infty}$ and module maps $\varphi_i :
F_i \raro F_{i-1}$, for all $i \geq 1$, such that the sequence of
modules
\[\cdots \longrightarrow F_m \stackrel{\varphi_m}{\longrightarrow}
F_{m-1} \longrightarrow \cdots \longrightarrow F_1
\stackrel{\varphi_{1}} \longrightarrow F_0 \stackrel{\varphi_{0}}
\longrightarrow M \longrightarrow 0,\]is exact where $F_0/
\mbox{ran} \varphi_1 = M$ and hence that $\varphi_0$ is a
projection. The above resolution is said to be a \textit{finite
resolution of length $l$}, for some $l \geq 0$, if $F_{l+1} = \{0\}$
and $F_i \neq \{0\}$ for $0 \leq i \leq l$.

A celebrated result in commutative algebra, namely, the Hilbert
Syzygy theorem, states that: Every finitely generated graded
$\mathbb{C}[\bm{z}]$-module has a finite graded free resolution of
length $l$ for some $l \leq n$ by finitely generated free modules.

It is also a question of general interest: given a free resolution
of a module over $\mathbb{C}[\bm{z}]$ when does the resolution stop.

\subsection{Dilations and free resolutions}

A contractive Hilbert module $\clh$ over $A(\mathbb{D})$ is said to
be \textit{completely non-unitary} (or c.n.u.) if there is no
non-zero reducing submodule $\cls \subseteq \clh$ such that
$M|_{\cls}$ is unitary.

Let $\clh$ be a contractive Hilbert module over $A(\mathbb{D})$.
Then the defect operators of $\clh$ are defined by $D_{\clh} =
(I_{\clh} - M^* M)^{\frac{1}{2}} \in \clb(\clh)$ and $D_{*\clh} =
(I_{\clh} - M M^*)^{\frac{1}{2}} \in \clb(\clh)$, and the defect
spaces by $\cld_{\clh} = \overline{\mbox{ran}} D_{\clh}$ and
$\cld_{*\clh} = \overline{\mbox{ran}} D_{*\clh}$. The
\textit{characteristic function} $\Theta_{\clh} \in
H^{\infty}_{\clb(\cld_{\clh}, \cld_{*\clh})} (\mathbb{D})$ is
defined by
\[\Theta_{\clh} (z) = [ - M + z D_{*\clh} (I_{\clh} - z M^*)^{-1}
D_{\clh}]|_{\cld_{\clh}}. \quad (z \in \mathbb{D})\] Define
$\Delta_{\clh}(t) = [ I_{\cld_{\clh}} - \Theta_{\clh}(e^{it})^*
\Theta_{\clh}(e^{it})]^{\frac{1}{2}} \in
\clb(L^2_{\cld_{\clh}}(\mathbb{T}))$ for $t \in [0, 1]$.
Then\[\clm_{\clh} = H^2_{\cld_{*\clh}}(\mathbb{D}) \oplus
\overline{\Delta_{\clh} L^2_{\cld_{\clh}}(\mathbb{T})},\]is a
contractive Hilbert module over $A(\mathbb{D})$. Then
\[\cls_{\clh} = \{M_{\Theta_{\clh}} f \oplus \Delta_{\clh} f : f \in
H^2_{\cld_{\clh}}(\mathbb{D})\} \subseteq \clm_{\clh},\]defines a
submodule of $\clm_{\clh}$. Also consider the quotient module
\[\mathcal{Q}_{\clh} = \clm_{\clh} \ominus \cls_{\clh}.\]
Here the module map $M_z \oplus M_{e^{it}}|_{\overline{\Delta_{\clh}
L^2_{\cld_{\clh}}(\mathbb{T})}}$ on $\clm_{\clh}$ is an isometry
where $M_z$ on $H^2_{\cld_{*\clh}}(\mathbb{D})$ is the pure part and
$M_{e^{it}}|_{\overline{\Delta_{\clh}
L^2_{\cld_{\clh}}(\mathbb{T})}}$ on $\overline{\Delta_{\clh}
L^2_{\cld_{\clh}}(\mathbb{T})}$ is the unitary part in the sense of
the Wold decomposition of isometries, Theorem \ref{Wold}.
Consequently,
\[ 0 \longrightarrow H^2_{\cld_{\clh}}(\mathbb{D})
\stackrel{\begin{bmatrix}\Theta_{\clh}
\\\Delta_{\clh}
\end{bmatrix}}{\longrightarrow} \clm_{\clh}
\,\;\stackrel{\Pi^*_{NF}}{\longrightarrow} \clq_{\clh}
\longrightarrow 0,\]where $\Pi_{NF}^*$ is the quotient (module) map.

\begin{Theorem}\label{Nagy-Foias-model}\textsf{(Sz.-Nagy and
Foias)} Let $\clh$ be a c.n.u. contractive Hilbert module over
$A(\mathbb{D})$. Then

(i) $\clh \cong \clq_{\clh}$.

(ii) $\clm_{\clh}$ is the minimal isometric dilation of $\clh$.
\end{Theorem}

Minimality of Sz.-Nagy-Foias isometric dilation, conclusion (ii) in
Theorem \ref{Nagy-Foias-model}, can be interpreted as a
factorization of dilation maps in the following sense:

\NI Let $\clh$ be a c.n.u. contractive Hilbert module over
$A(\mathbb{D})$ and $\Pi : \clh \raro \clk$ be an isometric dilation
of $\clh$ with isometry $V$ on $\clk$. Then there exists a unique
co-module isometry $\Phi \in \clb(\clm_{\clh}, \clk)$ such that
\[\Pi = \Phi \Pi_{NF},\]that is, the following diagram commutes:

 \setlength{\unitlength}{3mm}
 \begin{center}
 \begin{picture}(40,16)(0,0)
\put(15,3){$\clh$}\put(19,1.6){$\Pi$} \put(22.9,3){$\clk$} \put(22,
10){$\clm_{\clh}$} \put(22.8,9.2){ \vector(0,-1){5}} \put(15.8,
4.3){\vector(1,1){5.8}} \put(16.4,
3.4){\vector(1,0){6}}\put(16.5,8){$\Pi_{NF}$}\put(24.3,7){$\Phi$}
\end{picture}
\end{center}

As will be shown below, specializing to the case of $C_{\cdot 0}$
class and using localization technique one can recover the
characteristic function of a given $C_{\cdot 0}$-contractive Hilbert
module. Recall that a contractive Hilbert module $\clh$ over
$A(\mathbb{D})$ is said to be in $C_{\cdot 0}$ class if $M^{*k}
\raro 0$ in SOT as $k \raro \infty$. Submodules and quotient modules
of vector-valued Hardy modules are examples of Hilbert modules in
$C_{\cdot 0}$ class.

Let $\clh$ be a $C_{\cdot 0}$ contractive Hilbert module over
$A(\mathbb{D})$. Then there exists a Hilbert space $\cle_*$ such
that $\clh \cong \clq$ for some quotient module $\clq$ of
$H^2_{\cle_*}(\mathbb{D})$ (cf. Corollary \ref{pure-H2n}). Now by
Beurling-Lax-Halmos theorem, Theorem \ref{BLH-D}, there exists a
Hilbert space $\cle$ such that the submodule $\clq^{\perp} \cong
H^2_{\cle}(\mathbb{D})$. This yields the following short exact
sequence of modules:
\[0 \longrightarrow H^2_{\cle}(\mathbb{D})
\stackrel{X}{\longrightarrow} H^2_{\cle_*}(\mathbb{D})
\stackrel{\pi}{\longrightarrow} \clh \longrightarrow 0,\]where $X$
is isometric module map, and $\pi$ is co-isometric module map.
Localizing the isometric part of the short exact sequence,
$H^2_{\cle}(\mathbb{D}) \stackrel{X}{\longrightarrow}
H^2_{\cle_*}(\mathbb{D})$, at $z \in \mathbb{D}$ one gets
\[H^2_{\cle}(\mathbb{D})/ (A(\mathbb{D})_z \cdot
H^2_{\cle}(\mathbb{D})) \stackrel{X_z}{\longrightarrow}
H^2_{\cle_*}(\mathbb{D})/ (A(\mathbb{D})_z \cdot
H^2_{\cle_*}(\mathbb{D})).\]Identifying $H^2_{\cle}(\mathbb{D})/
(A(\mathbb{D})_z \cdot H^2_{\cle}(\mathbb{D}))$ with $\cle$ and
$H^2_{\cle_*}(\mathbb{D})/ (A(\mathbb{D})_z
H^2_{\cle_*}(\mathbb{D}))$ with $\cle_*$ one can recover the
characteristic function of $\clh$ as the map $z \mapsto X_z \in
\clb(\cle, \cle_*)$.

\subsection{Invariants} This subsection begins by proving a theorem, due to Sz.-Nagy and Foias
(\cite{NF}), on a complete unitary invariant of c.n.u. contractions.

\begin{Theorem}\label{coincidence}
Let $\clh_1$ and $\clh_2$ be c.n.u. contractive Hilbert modules over
$A(\mathbb{D})$. Then $\clh_1 \cong \clh_2$ if and only if
$\Theta_{\clh_1} \cong \Theta_{\clh_2}$.
\end{Theorem}

\NI \textsf{Proof.} Denote the module multiplication operator on
$\clh_1$ and $\clh_2$ by $M_1$ and $M_2$, respectively. Now let $u
M_1 = M_2 u$, for some unitary $u : \clh_1 \raro \clh_2$. Since $u
D_{* \clh_1} = D_{* \clh_2} u$ and $u D_{\clh_1} = D_{\clh_2} u$
\[ u|_{\cld_{\clh_1}} : \cld_{\clh_1} \raro \cld_{\clh_2} \quad \mbox{and}
\quad u|_{\cld_{* \clh_1}} : \cld_{* \clh_1} \raro \cld_{*
\clh_2},\]are unitary operators. A simple computation now reveals
that
\[u|_{\cld_{* \clh_1}} \Theta_{\clh_1}(z) = \Theta_{\clh_2}(z)
u|_{\cld_{\clh_1}},\]for all $z \in \mathbb{D}$, that is,
$\Theta_{\clh_1} \cong \Theta_{\clh_2}$.

\NI Conversely, given unitary operators $u \in \clb(\cld_{\clh_1},
\cld_{\clh_2})$ and $u_* \in \clb(\cld_{*\clh_1}, \cld_{*\clh_2})$
with the intertwining property $u_* \Theta_{\clh_1}(z) =
\Theta_{\clh_2}(z) u$ for all $z \in \mathbb{D}$,
\[\bm{u} = I_{H^2(\mathbb{D})} \otimes u|_{\cld_{\clh_1}} :
H^2_{\cld_{{\clh_1}}}(\mathbb{D}) \raro
H^2_{\cld_{{\clh_2}}}(\mathbb{D}),
\]and \[\bm{u_*} = I_{H^2(\mathbb{D})} \otimes u|_{\cld_{* \clh_1}}
: H^2_{\cld_{{* \clh_1}}}(\mathbb{D}) \raro H^2_{\cld_{{*
\clh_2}}}(\mathbb{D}),\] and \[\bm{\tau} = (I_{L^2(\mathbb{T})}
\otimes u)|_{\overline{\Delta_{\clh_1}
L^2_{\cld_{\clh_1}}(\mathbb{T})}} : \overline{\Delta_{\clh_1}
L^2_{\cld_{\clh_1}}(\mathbb{T})} \raro \overline{\Delta_{\clh_2}
L^2_{\cld_{\clh_2}}(\mathbb{T})},\] are module maps. Moreover,
\[\bm{u_*} M_{\Theta_{\clh_1}} = M_{\Theta_{\clh_2}} \bm{u}.\]
Consequently, one arrives at the following commutative diagram

\[\begin{CD}
0 @>>>H^2_{\cld_{\clh_1}}(\mathbb{D}) @>
\begin{bmatrix}M_{\Theta_{\clh_1}}\\\Delta_{\clh_1}\end{bmatrix}
>>H^2_{\cld_{*\clh_1}}(\mathbb{D}) @> \Pi^*_{NF, 1}>>\clq_{\clh_1} @ >>> 0\\
@. @V \bm{u} VV@V \bm{u}_* \oplus \tau VV @V VV\\
0 @>>>H^2_{\cld_{\clh_2}}(\mathbb{D}) @>>
\begin{bmatrix}M_{\Theta_{\clh_2}}\\\Delta_{\clh_2}\end{bmatrix} > H^2_{\cld_{*\clh_2}}(\mathbb{D}) @>> {\Pi^*_{NF, 2}}>\clq_{\clh_2} @
>>> 0
\end{CD}\]where the third vertical arrow is given by the unitary operator
\[\Pi_{NF, 2}^* (\bm{u_*} \oplus \bm{\tau}) \Pi_{NF, 1} : \clq_{\clh_1}
\raro \clq_{\clh_2}.\] To see this, first note that
\[\begin{split}(\bm{u_*} \oplus \bm{\tau})(\mbox{ran} \Pi_{NF, 1}) &
= (\bm{u_*} \oplus \bm{\tau})((\mbox{ker} \Pi_{NF, 1}^*))^{\perp} =
(\bm{u_*} \oplus \bm{\tau})((\mbox{ran}
\begin{bmatrix}M_{\Theta_{\clh_1}}\\\Delta_{\clh_1}\end{bmatrix})^{\perp})\\
& = [(\bm{u_*} \oplus \bm{\tau})(\mbox{ran}
\begin{bmatrix}M_{\Theta_{\clh_1}}\\\Delta_{\clh_1}\end{bmatrix})]^{\perp}
= [\mbox{ran}
\begin{bmatrix}M_{\Theta_{\clh_2}}\\\Delta_{\clh_2}\end{bmatrix}]^{\perp}\\&
= \mbox{ran} \Pi_{\clh_2}.
\end{split}\]
Moreover, the unitary operator \[(\bm{u_*} \oplus
\bm{\tau})|_{\mbox{ran} \Pi_{NF, 1}} : \mbox{ran} \Pi_{NF, 1} \raro
\mbox{ran} \Pi_{NF, 2},\]is a module map. This completes the proof.
\qed

\NI \textbf{Further results and comments:}

\begin{enumerate}
\item All results presented in this section can be found in the book
by Sz.-Nagy and Foias \cite{NF}. Here the Hilbert module point of
view is slightly different from the classical one.

\item Theorem \ref{coincidence} is due to Sz.-Nagy and Foias
\cite{NF}.

\item For non-commutative tuples of operators, Theorems \ref{Nagy-Foias-model} and \ref{coincidence} were
generalized by Popescu \cite{P1} and Ball and Vinnikov \cite{BV}
(see also \cite{GP99}, \cite{BES}, \cite{BTInd}, \cite{BTIEOT},
\cite{GP06}, \cite{GP11}, \cite{GP11a}, \cite{V}, \cite{B-Mem}, and
the references therein).

\item The notion of isometric dilation of contractions is closely related to the invariant subspace
problem (see \cite{ChPa}, \cite{Bez}, \cite{RR}). The reader is
referred to \cite{JS4}, \cite{JS3} for further recent developments
in this area.

\item There are many other directions to the model theory (both in single and several variables) that are
not presented in this survey. For instance, coordinate free approach
by Douglas, Vasyunin and Nikolski, and the de Branges-Rovnyak model
by de Branges, Rovnyak, Ball and Dritschel. We recommend the
monographs by Nikolski \cite{Nik} which is a comprehensive source of
these developments.

\item The paper by Ball and Kriete \cite{BKr} contains a remarkable
connection between the Sz.-Nagy and Foias functional model and the
de Branges-Rovnyak model on the unit disc.

\end{enumerate}

\section{Submodules}\label{S}

This section contains classical theory of isometries on Hilbert
spaces, invariant subspaces of $M_z$ on $H^2(\mathbb{D})$ and some
more advanced material on this subject.

Let $S$ be an isometry on a Hilbert space $\clh$, that is, $S^* S =
I_{\clh}$. A closed subspace $\clw \subseteq \clh$ is said to be
\textit{wandering subspace} for $S$ if $S^k \clw \perp S^l \clw$ for
all $k, l \in \mathbb{N}$ with $k \neq l$, or equivalently, if $S^k
\clw \perp \clw$ for all $k \geq 1$. An isometry $S$ on $\clh$ is
said to be \textit{shift} if
\[\clh = \mathop{\bigoplus}_{k \geq 0} S^k \clw,\]for some wandering
subspace $\clw$ for $S$. Equivalently, an isometry $S$ on $\clh$ is
shift if and only if (see Theorem \ref{Wold} below)
\[\mathop{\bigcap}_{k=0}^\infty S^k \clh = \{0\}.\]

For a shift $S$ on $\clh$ with a wandering subspace $\clw$ one sees
that
\[\clh \ominus S \clh = \mathop{\bigoplus}_{k \geq 0} S^k
\clw \ominus S(\mathop{\bigoplus}_{k \geq 0} S^k \clw) =
\mathop{\bigoplus}_{k \geq 0} S^k \clw \ominus \mathop{\bigoplus}_{m
\geq 1} S^k \clw = \clw.\]In other words, wandering subspace of a
shift is uniquely determined by $\clw = \clh \ominus S \clh$. The
dimension of the wandering subspace of a shift is called the
\textit{multiplicity} of the shift.

As for the examples, the only invariant of a shift operator is its
multiplicity, that is, the wandering subspace, up to unitary
equivalence.

\subsection{von Neumann and Wold decomposition}

One of the most important results in operator algebras, operator
theory and stochastic processes is the Wold decomposition theorem
\cite{Wo} (see also page 3 in \cite{NF}), which states that every
isometry on a Hilbert space is either a shift, or a unitary, or a
direct sum of shift and unitary.

\begin{Theorem}\label{Wold}
Let $S$ be an isometry on $\clh$. Then $\clh$ admits a unique
decomposition $\clh = \clh_s \oplus \clh_u$, where $\clh_s$ and
$\clh_u$ are $S$-reducing subspaces of $\clh$ and $S|_{\clh_s}$ is a
shift and $S|_{\clh_u}$ is unitary. Moreover,
\[\clh_s = \mathop{\bigoplus}_{k=0}^\infty S^k \clw \quad \quad
\mbox{and} \quad \quad \clh_u = \bigcap_{k=0}^\infty S^k
\clh,\]where $\clw = \mbox{ran}(I - S S^*)$ is the wandering
subspace for $S$.
\end{Theorem}
\NI\textsf{Proof.} Let $\clw = \mbox{ran}(I - S S^*)$ be the
wandering subspace for $S$ and $\clh_s =
\mathop{\bigoplus}_{k=0}^\infty V^k \clw$. Consequently, $\clh_s$ is
a $S$-reducing subspace of $\clh$ and that $S|_{\clh_s}$ is an
isometry. On the other hand, for all $k \geq 0$,
\[\begin{split} (S^k \clw)^\perp & = (S^k \mbox{ran}(I - S S^*))^\perp =
\mbox{ran} (I - S^k (I - S S^*) S^{*k}) \\ & = \mbox{ran}[ (I - S^k
S^{*k}) + S^{k+1} S^{* \,k+1}] = \mbox{ran}(I - S^k S^{*k}) \oplus
\mbox{ran} S^{k+1} \\ & = (S^k \clh)^\perp \oplus S^{k+1} \clh.
\end{split}\]Therefore \[\clh_u := \clh_s^\perp =
\mathop{\bigcap}_{k=0}^\infty S^k \clh.\] Uniqueness of the
decomposition readily follows from the uniqueness of the wandering
subspace $\clw$ for $S$. This completes the proof. \qed

\begin{Corollary}\label{Wold1}
Let $\clh$ be a Hilbert module over $\mathbb{C}[z]$. If the module
multiplication $M$ on $\clh$ is a shift then there exists a Hilbert
space $\clw$ and a module isometry $U$ from $H^2_{\clw}(\mathbb{D})$
onto $\clh$.
\end{Corollary}
\NI\textsf{Proof.} Let $\clw$ be the wandering subspace for $M$.
Define \[U : H^2_{\clw}(\mathbb{D}) \raro \clh =
\mathop{\bigoplus}_{k=0}^\infty M^k \clw,\]by $U(z^k f) = M^k f$ for
all $f \in \clw$ and $k \in \mathbb{N}$. One can check that this is
indeed the isometric module map from $H^2_{\clw}(\mathbb{D})$ onto
$\clh$. \qed

\subsection{Submodules of $H^2_{\cle}(\mathbb{D})$}
The purpose of this subsection is to show that a submodule of
$H^2_{\cle}(\mathbb{D})$ is uniquely determined (up to unitary
multipliers) by inner multipliers. The present methodology applies
the von Neumann-Wold decomposition theorem, to the submodules of the
Hardy module $H^2_{\cle}(\mathbb{D})$ (see page 239, Theorem 2.1 in
\cite{FF} and \cite{DNYJ}).

\begin{Theorem}\label{BLH-D}\textsf{(Beurling-Lax-Halmos Theorem)}
Let $\cls$ be a submodule of the Hardy module
$H^2_{\cle}(\mathbb{D})$. Then there exists a closed subspace $\clf
\subseteq \cle$ such that \[\cls \cong H^2_{\clf}(\mathbb{D}).\]In
particular, there exists an inner function $\Theta \in
H^\infty_{\cll(\clf, \cle)}(\mathbb{D})$ such that $M_\Theta :
H^2_\clf(\mathbb{D}) \raro H^2_\cle(\mathbb{D})$ is a module
isometry and $\cls = \Theta H^2_{\clf}(\mathbb{D})$. Moreover,
$\Theta$ is unique up to a unitary constant right factor, that is,
if $\cls = \tilde{\Theta} H^2_{\tilde{\clf}}(\mathbb{D})$ for some
Hilbert space $\tilde{\clf}$ and inner function $\tilde{\Theta} \in
H^\infty_{\clb(\tilde{\clf}, \cle)}(\mathbb{D})$, then $\Theta =
\tilde{\Theta} W$ where $W$ is a unitary operator in $\clb(\clf,
\tilde{\clf})$.
\end{Theorem}
\NI\textsf{Proof.} Let $\cls$ be a submodule of
$H^2_\cle(\mathbb{D})$. Then \[\bigcap_{l=0}^\infty (M_z|_{\cls})^l
\cls \subseteq \bigcap_{l=0}^\infty M_z^l H^2_{\cle}(\mathbb{D}) =
\{0\}.\] By Corollary \ref{Wold1} there exists an isometric module
map $U$ from $H^2_\clf(\mathbb{D})$ onto $\cls \subseteq
H^2_{\cle}(\mathbb{D})$. Consequently, $U = M_{\Theta}$ for some
inner function $\Theta \in H^\infty_{\cll(\clf, \cle)}(\mathbb{D})$.
\qed

In the particular case of the space $\cle = \mathbb{C}$, the above
result recovers Beurling's characterization of submodules of
$H^2(\mathbb{D})$.

\begin{Corollary}\label{Beur}\textsf{(Beurling)} Let $\cls$ be a
non-zero submodule of $H^2(\mathbb{D})$. Then $\cls = \theta
H^2(\mathbb{D})$ for some inner function $\theta \in
H^{\infty}(\mathbb{D})$.
\end{Corollary}

Moreover, one also has the following corollary:

\begin{Corollary}\label{Beur-2}
Let $\cls_1$ and $\cls_2$ be submodules of $H^2(\mathbb{D})$. Then
$\cls_1 \cong \cls_2$.
\end{Corollary}

The conclusion of Beurling's theorem, Corollary \ref{Beur}, fails if
$H^2(\mathbb{D})$ is replaced by the Bergman module
$L^2_a(\mathbb{D})$. However, a module theoretic interpretation of
Beurling-Lax-Halmos theorem states that: Let $\cls$ be a closed
subspace of the "free module" $H^2(\mathbb{D}) \otimes \cle (\cong
H^2_{\cle}(\mathbb{D}))$. Then $\cls$ is a submodule of
$H^2_{\cle}(\mathbb{D})$ if and only if $\cls$ is also "free" with
$\cls \ominus z \cls$ as a generating set. Moreover, in this case
$\mbox{dim} [\cls \ominus \cls] \leq \mbox{dim~} \cle$. In
particular, the wandering subspace $\cls \ominus z \cls$ is a
generating set of $\cls$.

Recall that a bounded linear operator $T$ on a Hilbert space $\clh$
is said to have the \textit{wandering subspace property} if $\clh$
is generated by the subspace $\clw_T : = \clh \ominus T \clh$, that
is,
\[\clh = [\clw_T] = \overline{\mbox{span}}\{T^m \clw_T : m \in
\mathbb{N}\}.\]In that case $\clw_T$ is said to be a wandering
subspace for $T$.

The following statements, due to Aleman, Richter and Sundberg
\cite{ARS-acta}, assert that the same conclusion hold also in the
Bergman module $L^2_a(\mathbb{D})$.

\begin{Theorem}\label{ARS-Berg}
Let $\cls$ be a submodule of $L^2_a(\mathbb{D})$. Then
\[\cls = \mathop{\bigvee}_{k=0}^\infty z^k (\cls \ominus z \cls).\]
\end{Theorem}

The same conclusion holds for the weighted Bergman space $L^2_{a,
\alpha}(\mathbb{D})$ with weight $\alpha = 3$ \cite{Shim} but for
$\alpha > 3$, the issue is more subtle (see \cite{HP}, \cite{MR}).

Another important consequence of the Beurling-Lax-Halmos theorem is
the characterization of cyclic submodules of
$H^2_{\cle}(\mathbb{D})$: Let $f$ be a non-zero vector in
$H^2_{\cle}(\mathbb{D})$. Then the cyclic submodule of
$H^2_{\cle}(\mathbb{D})$ generated by $f$ (and denoted by $[f]$) is
isomorphic to $H^2(\mathbb{D})$.

\NI There is no analog of the preceding result for the Bergman
module:

\begin{Theorem}\label{D-Bergman}
There does not exists any submodule $\cls$ of $L^2_a(\mathbb{D})$
such that $\cls \cong [1 \oplus z]$, the cyclic submodule of
$L^2_a(\mathbb{D}) \oplus L^2_a(\mathbb{D}) (\cong L^2_a(\mathbb{D})
\otimes \mathbb{C}^2)$ generated by $1 \oplus z$.
\end{Theorem}

\NI\textsf{Proof.} Let $\cls$ be a submodule of $L^2_a(\mathbb{D})$
and $U$ be a module isometric isomorphism from $[1 \oplus z]$ onto
$\cls$. Let \[U(1 \oplus z) = f,\] for some $f \in
L^2_a(\mathbb{D})$. Then the fact that the closed support of
Lebesgue measure on $\mathbb{D}$ is $\mathbb{D}$ implies that
\[|f(z)|^2 = 1 + |z|^2. \quad \quad (z \in \mathbb{D})\]By Taylor series
expansion of $f(z)$ one can show this is impossible for any
holomorphic function $f$ on $\mathbb{D}$. \qed

In the language of Hilbert modules, Beurling-Lax-Halmos theorem says
that the set of all non-zero submodules of
$H^2_{\cle_*}(\mathbb{D})$ are uniquely determined by the set of all
module isometric maps from $H^2_{\cle}(\mathbb{D})$ to
$H^2_{\cle_*}(\mathbb{D})$ where $\cle$ is a Hilbert space so that
$\mbox{dim~} \cle \leq \mbox{dim~} \cle_*$. On the other hand, a
module map $U : H^2_{\cle}(\mathbb{D}) \raro
H^2_{\cle_*}(\mathbb{D})$ is uniquely determined by a multiplier
$\Theta \in H^{\infty}_{\clb(\cle, \cle_*)}(\mathbb{D})$ and that
$\Theta$ is inner if and only if $U$ is isometry (cf. \cite{NF}).
Consequently, there exists a bijective correspondence, modulo the
unitary group, between the set of all non-zero submodules of
$H^2_{\cle_*}(\mathbb{D})$ and the set of all isometric module maps
from $H^2_{\cle}(\mathbb{D})$ to $H^2_{\cle_*}(\mathbb{D})$, where
$\cle \subseteq \cle_*$ and the set of all inner multipliers $\Theta
\in H^{\infty}_{\clb(\cle, \cle_*)}(\mathbb{D})$, where $\cle
\subseteq \cle_*$.

\subsection{Submodules of $H^2_n$}
This subsection will show how to extend the classification result of
submodules of $H^2_{\cle}(\mathbb{D})$, the Beurling-Lax-Halmos
theorem, to $H^2_n \otimes \cle$. This important generalization was
given by McCullough and Trent \cite{MT}.

Recall that the Drury-Arveson module $H^2_n \otimes \cle$ is a
reproducing kernel Hilbert module corresponding to the kernel\[(\z,
\w) \mapsto (1 - \sum_{i=1}^n z_i \bar{w}_i)^{-1} I_{\cle},\] for
all $\z, \w \in \mathbb{B}^n$ (see Section \ref{HMHF}). A multiplier
$\Theta \in \clm_{\clb(\cle, \cle_*)}(H^2_n)$ is said to be
\textit{inner} if $M_{\Theta}$ is a partial isometry in $\cll(H^2_n
\otimes \cle, H^2_n \otimes \cle_*)$.

\begin{Theorem}\label{MT}
Let $\cls (\neq \{0\})$ be a closed subspace of $H^2_n \otimes
\cle_*$. Then $\cls$ is a submodule of $H^2_n \otimes \cle_*$ if and
only if
\[\cls = \Theta(H^2_n \otimes \cle),\]for some inner multiplier
$\Theta \in \clm_{\clb(\cle, \cle_*)}(H^2_n)$.
\end{Theorem}

\NI\textsf{Proof.} Let $\cls$ be a submodule of $H^2_n \otimes
\cle_*$  and $R_i = M_{z_i}|_{\cls}$, $i = 1, \ldots, n$. Then
\[\sum_{i=1}^n R_{i} R_{i}^* = \sum_{i=1}^n P_{\cls} M_{z_i}P_{\cls}
M_{z_i}^* P_{\cls} \leq \sum_{i=1}^n P_{\cls} M_{z_i} M_{z_i}^*
P_{\cls},\]and consequently,
\[\begin{split} P_{\cls} - \sum_{i=1}^n R_{i} R_{i}^* & = P_{\cls} -
\sum_{i=1}^n P_{\cls} M_{z_i}P_{\cls} M_{z_i}^* P_{\cls} \geq
P_{\cls} - \sum_{i = 1}^n P_{\cls} M_{z_i} M_{z_i}^* P_{\cls} \\ & =
P_{\cls}(I_{H^2_n \otimes \cle_*} - \sum_{i=1}^n M_{z_i} M_{z_i}^*)
P_{\cls}.\end{split}\]Define $K : \mathbb{B}^n \otimes \mathbb{B}^n
\raro \cll(\cle_*)$, a positive definite kernel, by
\[\langle K(\bm{z}, \bm{w}) x_l, x_m \rangle = \langle
(P_{\cls} - \sum_{i = 1}^n P_{\cls} M_{z_i} P_{\cls} M_{z_i}^*
P_{\cls} ) (k_n(\cdot, \bm{w}) \otimes x_l) , k_n(\cdot, \bm{z})
\otimes x_m \rangle\]where $\{x_l\}$ is a basis of $\cle_*$. By
Kolmogorov theorem, there exists a Hilbert space $\cle$, a function
$\Theta \in \clo(\mathbb{B}^n, \clb(\cle, \cle_*))$ such that
\[K(\z, \w) = \Theta(\z) \Theta(\w)^*. \quad \quad (\bm{z}, \bm{w} \in
\mathbb{B}^n)\]On the other hand, since \[P_{\cls} M_{z_i} P_{\cls}
M_{z_i}^* P_{\cls} = M_{z_i} P_{\cls} M_{z_i}^*,\]for $i = 1,
\ldots, n$, we have
\[\begin{split}\langle (P_{\cls} - \sum_{i = 1}^n P_{\cls} M_{z_i} P_{\cls}
M_{z_i}^* P_{\cls} ) & (k_n(\cdot, \bm{w}) \otimes x_l) , k_n(\cdot,
\bm{z}) \otimes x_m \rangle \\ & = \langle (P_{\cls} - \sum_{i =
1}^n M_{z_i} P_{\cls} M_{z_i}^* ) (k_n(\cdot, \bm{w}) \otimes x_l) ,
k_n(\cdot, \bm{z}) \otimes x_m \rangle \\ & = k_n^{-1}(\z, \w)
\langle P_{\cls} (k_n(\cdot, \bm{w}) \otimes x_l) , k_n(\cdot,
\bm{z}) \otimes x_m \rangle.
\end{split}\]Thus \[\langle k_n(\z, \w) K(\z, \w) x_l, x_m\rangle = \langle P_{\cls} (k_n(\cdot, \bm{w}) \otimes x_l) , k_n(\cdot,
\bm{z}) \otimes x_m \rangle.\]This implies that
\[(\z, \w) \mapsto \big(I_{\cle_*} - K(\z, \w)\big) k_n(\z, \w) = \big(I_{\cle_*} - \Theta(\z) \Theta(\w)^*\big) k_n(\z, \w)\]is
a $\clb(\cle_*)$-valued positive definite kernel, from which it
follows that $\Theta$ is a multiplier, that is, $\Theta \in
\clm_{\clb(\cle, \cle_*)}(H^2_n)$. Finally,
\[\begin{split}\langle M_{\Theta} M_{\Theta}^*
(k_n(\cdot, \bm{w}) \otimes x_l), k_n(\cdot, \bm{z}) \otimes x_m
\rangle & =  \langle k_n(\bm{z}, \bm{w}) \Theta(\bm{z})
\Theta(\bm{w})^* x_l, x_m\rangle\\ &  = \langle k_n(\bm{z}, \bm{w})
K(\bm{z}, \bm{w}) x_l, x_m\rangle \\
& = \langle P_{\cls} (k_n(\cdot, \bm{w}) \otimes x_l), k_n(\cdot,
\bm{z}) \otimes x_m\rangle,\end{split}\]and hence $P_{\cls} =
M_{\Theta} M_{\Theta}^*$ and that $M_{\Theta}$ a partial isometry.
This completes the proof. \qed

In \cite{GRS}, Green, Richter and Sundberg prove that for almost
every $\bm{\zeta} \in \partial \mathbb{B}^n$ the nontangential limit
$\Theta(\bm{\zeta})$ of the inner multiplier $\Theta$ is a partial
isometry. Moreover, the rank of $\Theta(\bm{\zeta})$ is equal to a
constant almost everywhere.

\subsection{Solution to a Toeplitz operator equation}
This subsection contains an application of Hilbert module approach
to a problem concerning the classical analytic Toeplitz operators.
This Toeplitz operator equation problem can be formulated in a more
general framework.

Let $\cls = M_{\Theta}H^2_{\cle}(\mathbb{D})$ be a $M_z$-invariant
subspace of $H^2_{\cle_*}(\mathbb{D})$ for some inner multiplier
$\Theta \in H^{\infty}_{\clb(\cle, \cle_*)}(\mathbb{D})$. Moreover,
let $\cls$ be invariant under $M_{\Phi}$ for some $\Phi \in
\clb(\cle_*)$. Then
\[ \Phi \,\Theta = \Theta \Psi,\]for some unique $\Psi \in H^{\infty}_{\clb(\cle)}(\mathbb{D})$.

\NI \textsf{Problem:} Determine $\Psi$, that is, find a
representation of the unique multiplier $\Psi$. If $\Phi$ is a
polynomial, then under what conditions will $\Psi$ be a polynomial,
or a polynomial of the same degree as $\Phi$?

More precisely, given $\Theta$ and $\Phi$ as above, one seeks a
(unique) solution $X \in H^{\infty}_{\clb(\cle)}(\mathbb{D})$ to the
Toeplitz equation $\Theta X = \Phi \,\Theta$.

\NI This problem appears to be difficult because there are
infinitely many obstructions (rather, equations, if one expands
$\Theta$ and $\Phi$ in power series). Thus a priori the answer is
not expected to be tractable in general. However, it turns out that
if $\Phi(z) = A + A^* z$, then $\Psi = B + B^* z$ for some unique
$B$. The proof is a straightforward application of methods
introduced by Agler and Young in \cite{AY}. However, the intuitive
idea behind this "guess" is that, $\Phi$ turns
$H^2_{\cle_*}(\mathbb{D})$ into a natural Hilbert module over
$\mathbb{C}[z_1, z_2]$ (see Corollary \ref{BLH-mod}).

It is now time to proceed to the particular framework for the
Toeplitz operator equation problem. Let \[\Gamma = \{(z_1 + z_2, z_1
z_2): |z_1|, |z_2| \leq 1\} \subseteq \mathbb{C}^2,\] be the
symmetrized bidisc. A Hilbert module $\clh$ over $\mathbb{C}[z_1,
z_2]$ is said to be \textit{$\Gamma$-normal Hilbert module} if $M_1$
and $M_2$ are normal operators and $\sigma_{Tay}(M_1, M_2)$, the
Taylor spectrum of $(M_1, M_2)$ (see Section \ref{SFHM}), is
contained in the distinguished boundary of $\Gamma$. A Hilbert
module $\clh$ over $\mathbb{C}[z_1, z_2]$ is said to be
\textit{$\Gamma$-isometric Hilbert module} if $\clh$ is a submodule
of a $\Gamma$-normal Hilbert module. A $\Gamma$-isometric Hilbert
module $\clh$ is \textit{pure} if $M_2$ is a shift operator.

Let $\cle_*$ be a Hilbert space and $A \in \clb(\cle_*)$ with
$w(A)$, the numerical radius of $A$, not greater than one. By
$[H^2_{\cle_*}(\mathbb{D})]_A$ we denote the Hilbert module
$H^2_{\cle_*}(\mathbb{D})$ with \[\mathbb{C}[z_1, z_2] \times
H^2_{\cle_*}(\mathbb{D}) \raro H^2_{\cle_*}(\mathbb{D}), \quad \quad
(p(z_1, z_2), h) \mapsto p(A + A^* M_z, M_z)h.\]

\NI The following theorem is due to Agler and Young (see \cite{AY}).

\begin{Theorem}\label{gamma-isometry}
Let $\clh$ be a Hilbert module over $\mathbb{C}[z_1, z_2]$. Then
$\clh$ is a pure $\Gamma$-isometric Hilbert module if and only if
$\clh \cong [H^2_{\cle_*}(\mathbb{D})]_A$ for some Hilbert space
$\cle_*$, $A \in \clb(\cle_*)$ and $w(A) \leq 1$.
\end{Theorem}
Given a Hilbert space $\cle_*$ and  $A \in \clb(\cle_*)$ with $w(A)
\leq 1$, the Hilbert module $[H^2_{\cle_*}(\mathbb{D})]_A$ is called
a \textit{$\Gamma$-isometric Hardy module} with symbol $A$.

Now let $\cls$ be a non-zero submodule of
$[H^2_{\cle}(\mathbb{D})]_A$. Then in particular, by the
Beurling-Lax-Halmos theorem, Theorem \ref{BLH-D}, we have  \[\cls =
\Theta H^2_{\cle}(\mathbb{D}),\] for some Hilbert space $\cle$ and
inner multiplier $\Theta \in H^\infty_{\clb(\cle,
\cle_*)}(\mathbb{D})$.

Now everything is in place to state and prove the main result of
this subsection.

\begin{Theorem}\label{BLH-gamma}
Let $\cls \neq \{0\}$ be a closed subspace of
$H^2_{\cle_*}(\mathbb{D})$ and $A \in \clb(\cle_*)$ with $w(A) \leq
1$. Then $\cls$ is a submodule of $[H^2_{\cle_*}(\mathbb{D})]_A$ if
and only
\[(A + A^* M_z) M_{\Theta} = M_{\Theta} (B + B^* M_z),\]for some unique $B \in
\clb(\cle)$ (up to unitary equivalence) with $w(B) \leq 1$ where
$\Theta \in H^{\infty}_{\clb(\cle, \cle_*)}(\mathbb{D})$ is the
Beurling-Lax-Halmos representation of $\cls$.
\end{Theorem}

\NI\textsf{Proof.} Assume that $\cls$ be a non-zero submodule of
$[H^2_{\cle_*}(\mathbb{D})]_A$ and $\cls = M_{\Theta}
H^2_{\cle}(\mathbb{D})$ be the Beurling-Lax-Halmos representation of
$\cls$ where $\Theta \in H^{\infty}_{\clb(\cle,
\cle_*)}(\mathbb{D})$ is an inner multiplier and $\cle$ is an
auxiliary Hilbert space. Also
\[(A + A^* M_z) (M_{\Theta} H^2_{\cle}(\mathbb{D})) \subseteq
M_{\Theta} H^2_{\cle}(\mathbb{D}),\]implies that $(A + A^* M_z)
M_{\Theta} = M_{\Theta} M_{\Psi}$ for some unique $\Psi \in
H^{\infty}_{\clb(\cle)}(\mathbb{D})$. Therefore,\[M_{\Theta}^* (A +
A^* M_z) M_{\Theta} = M_{\Psi}.\]Multiplying both sides by $M_z^*$,
one arrives at
\[M_z^* M_{\Theta}^* (A + A^* M_z) M_{\Theta} = M_z^*
M_{\Psi}.\]Then $M_{\Theta}^* (A M_z^* + A^*) M_{\Theta} = M^*_z
M_{\Psi}$ and hence, $M_z^* M_{\Psi} = M_{\Psi}^*$, or equivalently,
$M_{\Psi} = M_{\Psi}^* M_z$. Since $\|M_{\Psi}\| \leq 2$, it follows
that $(M_{\Psi}, M_z)$ is a $\Gamma$-isometry. By Theorem
\ref{gamma-isometry}, it follows that
\[M_{\Psi} = B + B^* M_z,\]for some $B \in \clb(\cle)$ and $w(B) \leq
1$, and uniqueness of $B$ follows from that of $\Psi$.

\NI The converse part is trivial, and the proof is complete. \qed

One of the important applications of the above theorem is the
following result concerning Toeplitz operators with analytic
polynomial symbols of the form $A + A^* z$.

\begin{Theorem}\label{BLH-mult}
Let $\cls = M_{\Theta} H^2_{\cle}(\mathbb{D}) \subseteq
H^2_{\cle_*}(\mathbb{D})$ be a non-zero $M_z$-invariant subspace of
$H^2_{\cle_*}(\mathbb{D})$ and $A \in \clb(\cle_*)$. Then $\cls$ is
invariant under the Toeplitz operator with analytic polynomial
symbol $A + A^*z$ if and only if there exists a unique operator $B
\in \clb(\cle)$ such that
\[(A + A^* z) \Theta = \Theta (B + B^* z).\]
\end{Theorem}
The following result relates Theorem \ref{BLH-gamma} to module maps
of $\Gamma$-isometric Hardy modules.

\begin{Corollary}\label{BLH-mod}
Let $\cls \neq \{0\}$ be a closed subspace of
$H^2_{\cle_*}(\mathbb{D})$. Then $\cls$ is a submodule of the
$\Gamma$-isometric Hardy module $[H^2_{\cle_*}(\mathbb{D})]_A$ with
symbol $A$ if and only if there exists a $\Gamma$-isometric Hardy
module $[H^2_{\cle}(\mathbb{D})]_B$ with a unique symbol $B \in
\clb(\cle)$ and an isometric module map
\[U : [H^2_{\cle}(\mathbb{D})]_B \longrightarrow
[H^2_{\cle_*}(\mathbb{D})]_A,\]such that $\cls = U
H^2_{\cle}(\mathbb{D})$.
\end{Corollary}
Another application of Theorem \ref{BLH-gamma} concerns unitary
equivalence of $\Gamma$-isometric Hardy module submodules.

\begin{Corollary}
A non-zero submodule of a $\Gamma$-isometric Hardy module is
isometrically isomorphic with a $\Gamma$-isometric Hardy module.
\end{Corollary}

\NI \textbf{Further results and comments:}
\begin{enumerate}

\item The classification result of invariant subspaces, Corollary \ref{Beur},
is due to Beurling \cite{B}. The Beurling-Lax-Halmos theorem was
obtained by Lax \cite{L} and Halmos \cite{H} as a generalization of
Beurling's theorem (see \cite{NF}). See also the generalization by
Ball and Helton in \cite{BH}.

\NI The simple proof of the Beurling-Lax-Halmos theorem presented
here requires the von Neumann-Wold decomposition theorem which
appeared about two decades earlier than Beurling's classification
result on invariant subspaces of $H^2(\mathbb{D})$.

\item Let $\cls \neq \{0\}$ be a submodule of $H^2(\mathbb{D})$. Then the
wandering subspace of $\cls$, $\cls \ominus z \cls$, has dimension
one. However, in contrast with the Hardy module $H^2(\mathbb{D})$,
the dimension of the generating subspace $\cls \ominus z \cls$ of a
submodule $\cls$ of the Bergman module $L^2_a(\mathbb{D})$ could be
any number in the range $1, 2, \ldots$ including $\infty$. This
follows from the dilation theory developed by Apostol, Bercovici,
Foias and Pearcy (see \cite{ABFP}).

\item Beurling type theorem for the Bergman space,
Theorem \ref{ARS-Berg}, is due to Aleman, Richter and Sundberg. This
result was further generalized by Shimorin \cite{Shim} in the
context of operators close to isometries. His results include the
Dirichlet space on the unit disc. A several variables analogue of
the wandering subspace problem for the Bergman space over
$\mathbb{D}^n$ is proposed in \cite{CDSS}.

\item See \cite{III-Arc} for a simple and ingenious proof of the
Aleman-Richter-Sundberg theorem concerning invariant subspaces of
the Bergman space.

\item The proof of Theorem \ref{Wold} is from \cite{S-W}. It is slightly
simpler than the one in \cite{NF} and \cite{FF}. Theorem
\ref{BLH-gamma} and Corollary \ref{BLH-mod} are due to the author.
Theorem \ref{D-Bergman} is due to Douglas (\cite{DNYJ}).

\item Theorem \ref{MT} is due to McCollough and Trent \cite{MT}. For more related results in
one variable, see the article by Jury \cite{MJ}. See \cite{JS4} for
a new approach to Theorem \ref{MT}.

\item One possible approach to solve the problem mentioned in the last subsection is
to consider first the finite dimension case, that is, $\cle_* =
\mathbb{C}^k$ for $k>1$.

\item
Let $\cls \neq \{0\}$ be a closed subspace of
$H^2_{H^2(\mathbb{D}^{n-1})}(\mathbb{D})$. By Beurling-Lax-Halmos
theorem, that $\cls$ is a submodule of
$H^2_{H^2(\mathbb{D}^{n-1})}(\mathbb{D})$ if and only if $\cls =
\Theta H^2_{\cle_*}(\mathbb{D})$, for some closed subspace $\cle_*
\subseteq H^2(\mathbb{D}^{n-1})$ and inner function $\Theta \in
H^\infty_{\cll(\cle_*, H^2(\mathbb{D}^{n-1}))}(\mathbb{D})$. Here
one is naturally led to formulate the following problem.

\NI\textsf{Problem:} For which closed subspace $\cle_* \subseteq
H^2(\mathbb{D}^{n-1})$ and inner function $\Theta \in
H^\infty_{\clb(\cle_*, H^2(\mathbb{D}^{n-1}))}(\mathbb{D})$ the
submodule $\Theta H^2_{\cle_*}(\mathbb{D})$ of
$H^2_{H^2(\mathbb{D}^{n-1})}(\mathbb{D})$, realized as a subspace of
$H^2(\mathbb{D}^n)$, is a submodule of $H^2(\mathbb{D}^n)$?

\NI This problem is hard to tackle in general. However, see
\cite{JS1} for some partial results.

\item The Beurling-Lax-Halmos theorem for submodules of
vector-valued Hardy modules can be restated by saying that the
non-trivial submodules of $H^2_{\cle_*}(\mathbb{D})$ are the images
of vector-valued Hardy modules under partially isometric module maps
(see \cite{RR}). This classification result for $C_{\cdot
0}$-contractive Hilbert modules over $A(\mathbb{D})$ have also been
studied (see \cite{JS3}).

\begin{Theorem}
Let $\clh$ be a $C_{\cdot 0}$-contractive Hilbert module over
$A(\mathbb{D})$ and $\cls$ be a non-trivial closed subspace of
$\clh$. Then $\cls$ is a submodule of $\clh$ if and only if there
exists a Hilbert space $\cle$ and a partially isometric module map
$\Pi : H^2_{\cle}(\mathbb{D}) \raro \clh$ such that
\[\cls = \mbox{ran~} \Pi,\]or equivalently, \[P_{\cls} = \Pi
\Pi^*.\]
\end{Theorem}
\NI An analogous assertion is true also for Hilbert modules over
$\mathbb{C}[\z]$ (see \cite{JS4}).

\item Let $\cle$ and $\cle_*$ be Hilbert spaces and $m \in \mathbb{N}$. Let
$\Theta \in \clm(H^2_{\cle}(\mathbb{D}), L^2_{a, m}(\mathbb{D})
\otimes \cle_*)$ be a partially isometric multiplier. It follows
easily from the definition of multipliers that $\Theta
H^2_{\cle}(\mathbb{D})$ is a submodule of $L^2_{a, m}(\mathbb{D})
\otimes \cle_*$. The following converse was proved by Ball and
Bolotnikov in \cite{BB1} (see also \cite{BB2} and Olofsson
\cite{O}).

\begin{Theorem}\label{BBthm}
Let $\cls$ be a non-trivial submodule of the vector-valued weighted
Bergman module $L^2_{a, m}(\mathbb{D}) \otimes \cle_*$. Then there
exists a Hilbert space $\cle$ and partially isometric multiplier
$\Theta \in \clm(H^2_{\cle}(\mathbb{D}), L^2_{a, m}(\mathbb{D})
\otimes \cle_*)$ such that
\[\cls = \Theta H^2_{\cle}(\mathbb{D}).\]
\end{Theorem}
\NI Another representation for $\cls$, a submodule of $L^2_{a,
m}(\mathbb{D}) \otimes \cle_*$, is based on the observation that for
any such $\cls$, the subspace $z^k \cls \ominus z^{k+1} \cls$ can be
always represented as $z^k \Theta_{k} \clu_k$ for an appropriate
subspace $\clu_k$ and an $L^2_{a, m}(\mathbb{D}) \otimes
\cle_*$-inner function $z^k \Theta_{k}$, $k \geq 0$. This
observation leads to the orthogonal representation: \[\cls =
\oplus_{k \geq 0} (z^k \cls \ominus z^{k+1} \cls) = \oplus_{k \geq
0} z^k \Theta_{k} \clu_k,\]of $\cls$ in terms of a Bergman-inner
family $\{\Theta_k\}_{k \geq 0}$ (see \cite{BB1} and \cite{BB2} for
more details).

\NI More recently, Theorem \ref{BBthm} has been extended by the
author \cite{JS4}, \cite{JS3} to the case of reproducing kernel
Hilbert modules.

\end{enumerate}

\section{Unitarily equivalent submodules}\label{UES}

Let $\clh \subseteq \clo(\mathbb{D}, \mathbb{C})$ be a reproducing
kernel Hilbert module and $\cls_1$ and $\cls_2$ be two non-zero
submodules of $\clh$.

\begin{enumerate}
\item If $\clh = H^2(\mathbb{D})$, then $\cls_1 \cong \cls_2$ (see Corollary
\ref{Beur-2}).

\item If $\clh = L^2_a(\mathbb{D})$ and $\cls_1 \cong \cls_2$, then
$\cls_1 = \cls_2$ (see \cite{Rich} or Corollary 8.5 in
\cite{S-HM}).
\end{enumerate}
Therefore, on the one hand every non-zero submodule is isometrically
isomorphic to the module itself while on the other hand, no proper
submodule is.

\NI Now let $n> 1$. For submodules of $H^2(\mathbb{D}^n)$ over
$A(\mathbb{D}^n)$, some are unitarily equivalent to
$H^2(\mathbb{D}^n)$ and some are not (cf. \cite{R}, \cite{M-PAMS},
\cite{SSW}). For the Hardy module $H^2(\partial\mathbb{B}^n)$, the
existence of inner functions on $\mathbb{B}^n$ \cite{Alek}
established the existence of proper submodules of
$H^2(\partial\mathbb{B}^n)$ that are unitarily equivalent to
$H^2(\mathbb{B}^n)$.

These observations raise a number of interesting questions
concerning Hilbert modules with unitarily equivalent submodules. The
purpose of this section is to investigate and classify a class of
Hilbert modules with proper submodules unitarily equivalent to the
original.

\subsection{Isometric module maps}
This subsection begins with a simple observation concerning
unitarily equivalent submodules of Hilbert modules. Let $\clh$ be a
Hilbert module over $A(\Omega)$ and $\cls$ be a non-trivial
submodule of $\clh$. Then $\cls$ is unitarily equivalent to $\clh$
if and only if $\cls = U \clh$ for some isometric module map $U$ on
$\clh$.

Now let $U\clh$ be a submodule of $\clh$ for some isometric module
map $U$. Then $U \clh$ is said to be \emph{pure} unitarily
equivalent submodule of $\clh$ if\[\bigcap\limits^\infty_{k\ge 0}
U^k \mathcal{H}=\{0\}.\]

\begin{Proposition}\label{pro1}
Let $\mathcal{H}$ be a Hilbert module over $A(\Omega)$ for which
there exists an isometric module map $U$ satisfying
$\bigcap\limits^\infty_{k=0} U^k\mathcal{H} = (0)$. Then there
exists an isomorphism $\Psi \colon \
H^2_{\mathcal{W}}(\mathbb{D})\to \clh$ with $\mathcal{W} = \clh
\ominus U \clh$ and a commuting $n$-tuple of functions
$\{\varphi_i\}$ in $H^\infty_{\clb(\mathcal{W})}(\mathbb{D})$ so
that $U = \Psi M_z\Psi^*$ and $M_{i} = \Psi M_{\varphi_i}\Psi^*$ for
$i=1,2,\ldots, n$.
\end{Proposition}
\NI\textsf{Proof.} By Corollary \ref{Wold1}, there is a canonical
isomorphism $\Psi\colon \ H^2_{\mathcal{W}} (\mathbb{D})\to \clh$
such that $\Psi T_z = U\Psi$ where $\clw = \clh \ominus U \clh$.
Further, $X_i = \Psi^*M_{i}\Psi$ is an operator on
$H^2_{\mathcal{E}}(\mathbb{D})$ which commutes with $T_z$. Hence,
there exists a function $\varphi_i$ in
$H^\infty_{\mathcal{L}(\mathcal{W})}(\mathbb{D})$ such that $X_i =
M_{\varphi_i}$. Moreover, since the $\{M_{i}\}$ commute, so do the
$\{X_i\}$ and hence the functions $\{\varphi_i\}$ commute pointwise
a.e.\ on $\mathbb{T}$. \qed

\subsection{Hilbert-Samuel Polynomial}

A Hilbert module $\clh$ over $\mathbb{C}[z]$ is said to be
semi-Fredholm at $\w \in \mathbb{C}^n$ if \[\mbox{dim} [\clh/I_{\w}
\cdot \clh] < \infty.\]In particular, note that $\clh$ semi-Fredholm
at $\w$ implies that $I_{\w} \cdot \clh$ is a closed submodule of
$\clh$ and\[\mbox{dim} [I^k_{\w}\cdot \clh/ I^{k+1}_{\w} \cdot \clh]
< \infty,\]for all $k \in \mathbb{N}$. In this case the direct sum
\[\mbox{gr}(\clh) : = \mathop{\bigoplus}_{k \geq 0} I^k_{\w}\cdot
\clh/ I^{k+1}_{\w} \cdot \clh,\]can be turned into a graded finitely
generated $\mathbb{C}[\z]$-module. It is a fundamental result of
commutative algebra that to any such module there is a polynomial
$h^{\w}_{\clh} \in \mathbb{Q}[x]$ of degree not greater than $n$,
the \textit{Hilbert-Samuel polynomial}, with\[h^{\w}_{\clh}(k) =
\mbox{dim} [I^{k}_{\w} \cdot \clh/ I^{k+1}_{\w} \cdot \clh],\]for
all $k \geq N_{\clh}$ for some positive integer $N_{\clh}$ (see
\cite{D-Y}).

In some cases it is possible to calculate the Hilbert-Samuel
polynomial for a Hilbert module directly. For example (see
\cite{F-05}), let $\Omega$ be a Reinhardt domain in $\mathbb{C}^n$
and $\clh \subseteq \clo(\Omega, \mathbb{C})$ be a reproducing
kernel Hilbert module. Let $\cls$ be a singly generated submodule of
$\clh$ and $\w \in \Omega$. Then
\[h_{\cls}^{\w} (k) = \binom{n+k - 1}n.\]In general, it is difficult
to compute the Hilbert-Samuel polynomial directly.

The following result demonstrates that the Hilbert-Samuel polynomial
does not depend on the choice of a quasi-free Hilbert module.
\begin{Theorem}\label{lemma-DS}
If $\mathcal{R}$ and $\widetilde{\mathcal{R}}$ be a pair of rank $m$
quasi-free Hilbert modules over $A(\Omega)$. If both $\clr$ and
$\widetilde{\mathcal{R}}$ are semi-Fredholm at $\w \in \Omega$ then
$h_{\mathcal{R}}^{\pmb{\omega}} \equiv
h_{\widetilde{\mathcal{R}}}^{\pmb{\omega}}$.
\end{Theorem}

\NI\textsf{Proof.}  Consider rank $m$ quasi-free Hilbert modules
$\mathcal{R}$ and $\widetilde{\mathcal{R}}$ over $A(\Omega)$ with
$1\le m <\infty$. Following Lemma $1$ in \cite{DM2}, construct the
rank $m$ quasi-free Hilbert module $\Delta$, which is the graph of a
closed densely defined module map from $\mathcal{R}$ to
$\widetilde{\mathcal{R}}$ obtained as the closure of the set
$\{\varphi f_i\oplus \varphi g_i\colon\ \varphi\in A(\Omega)\}$,
where $\{f_i\}_{i=1}^m$ and $\{g_i\}_{i=1}^m$ are generators for
$\mathcal{R}$ and $\widetilde{\mathcal{R}}$, respectively. Then the
module map $X\colon \ \Delta\to \mathcal{R}$ defined by $f_i\oplus
g_i \to f_i$ is bounded, one-to-one and has dense range. Note that
for fixed $\w_0$ in $\Omega$, \[X^*(I_{\w_0}\cdot \mathcal{R})^\bot
\subset (I_{\w_0}\cdot \Delta)^\bot.\]Since the rank of $\Delta$ is
also $k$, this map is an isomorphism. Let $\{\gamma_i(\w_0)\}$ be
anti-holomorphic functions from a neighborhood $\Omega_0$ of $\w_0$
to $\mathcal{R}$ such that $\{\gamma_i(\w)\}$ spans $(I_{\w} \cdot
\mathcal{R})^\bot$ for $\w \in \Omega_0$. Then
\[\{\frac{\partial^{\bm{k}}}{\partial
{z}^{\bm{k}}}\gamma_i(\w)\}_{|\bm{k}|< k},\] forms a basis for
$(I^k_{\w}\cdot \mathcal{R})^\bot$ for $k=0,1,2,\ldots$, using the
same argument as in Section 4 in \cite{CS} and Section 4 in
\cite{DMV}. Similarly, since $\{X^*\gamma_i(\w)\}$ is a basis for
$(I_{\w}\cdot \Delta)^\bot$, it follows that $X^*$ takes
$(I^k_{\w}\cdot \mathcal{R})^\bot$ onto $(I^k_{\w}\cdot
\Delta)^\bot$ for $k=0,1,2,\ldots$~. Therefore \[\dim(I^k_{\w}\cdot
\mathcal{R})^\bot = \dim(I^k_{\w}\cdot\Delta)^\bot,\]for all $k$.
Hence
\[h_{\mathcal{R}}^{\w} = h_{\Delta}^{\w}. \quad \quad (\w \in \Omega)\]
The result now follows by interchanging the roles of $\mathcal{R}$
and $\widetilde{\mathcal{R}}$. \qed

In particular, one can calculate the Hilbert-Samuel polynomial by
considering only the Bergman module over $A(\Omega)$ since
\[h_{\mathcal{R}\otimes \mathbb{C}^k}^{\w}\equiv k
h_{\mathcal{R}}^{\w},\]for all finite integers $k$. To accomplish
that one can reduce to the case of a ball as follows.

\begin{Theorem}\label{pro2}
If $\mathcal{R}$ is a quasi-free Hilbert module over $A(\Omega)$ for
$\Omega\subset \mathbb{C}^n$ which is semi-Fredholm for $\w$ in a
neighborhood of $\w_0$ in $\Omega$ with constant codimension, then
$h^{\w_0}_{\mathcal{R}}$ has degree  $n$.
\end{Theorem}
\NI\textsf{Proof.} Let $B_\varepsilon(\w_0)$ be a ball with radius
$\varepsilon$ centered at $\w_0$, whose closure is contained in
$\Omega$. An easy argument shows that the map $X\colon \
L^2_a(\Omega)\to L^2_a(B_\varepsilon(\w_0))$ defined by \[Xf\equiv
f|_{B_\varepsilon(\w_0)}, \quad \quad  (f \in L^2_a(\Omega))\]is
bounded, one-to-one and has dense range. Moreover, by a similar
argument to the one used in Theorem \ref{lemma-DS} for $\w \in
B_\varepsilon(\w_0)$, it follows that \[h_{L^2_a(\Omega)}^{\w}
\equiv h_{B_\varepsilon(\w)}^{\w}.\] The proof is completed by
considering the Hilbert--Samuel polynomials at $\w_0$ of the Bergman
module for the ball $B_\varepsilon(\w_0)$ for some $\varepsilon>0$
which is centered   at $\w_0$. This calculation reduces to that of
the module $\mathbb{C}[\z]$ over the algebra $\mathbb{C}[\pmb{z}]$
since the monomials in $L^2_a(B_\varepsilon(\w_0))$ are orthogonal.
Hence\[h_{L^2_a(B_\varepsilon(\w_0))}^{\w_0}(k) = \binom{n+k-1}n.\]
This completes the proof. \qed

\subsection{On complex dimension}

The purpose of this subsection is to show that the complex dimension
of the domain $\Omega$ is one, that is $n = 1$, whenever $\clh$ is
quasi-free, semi-Fredholm and $\dim \mathcal{H}/
U\mathcal{H}<\infty$.

The following result relates pure isometrically isomorphic
submodules of finite codimension and linear Hilbert-Samuel
polynomials.

\begin{Theorem}\label{pro3}
If $\mathcal{H}$ is semi-Fredholm at $\w_0$ in $\Omega$ and
$\mathcal{S}$ is a pure isometrically isomorphic submodule of
$\mathcal{H}$ having finite codimension in $\mathcal{H}$, then
$h^{\w_0}_{\mathcal{H}}$ has degree at most one.
\end{Theorem}

\NI\textbf{Proof.} As in the proof of Proposition \ref{pro1}, the
existence of $\mathcal{S}$ in $\mathcal{H}$ yields a module
isomorphism $\Psi$ of $\mathcal{H}$ with
$H^2_{\mathcal{W}}(\mathbb{D})$ for $\mathcal{W} =
\mathcal{H}\ominus \mathcal{S}$. Assume that $\w_0 = 0$ for
simplicity and note that the assumption that $\mathcal{H}$ is
semi-Fredholm at $\w_0 = 0$ implies that
\[M_{z_1}\cdot \clh +\cdots+ M_{z_n}\cdot \clh,\] has finite
codimension in $\mathcal{H}$. Hence \[\widetilde{\mathcal{S}} =
M_{\varphi_1}\cdot H^2_{\mathcal{W}}(\mathbb{D}) +\cdots+
M_{\varphi_n}\cdot H^2_{\mathcal{W}}(\mathbb{D}),\] has finite
codimension in $H^2_{\mathcal{W}}(\mathbb{D})$, where $M_{z_i} =
\Psi M_{\varphi_i} \Psi^*$. Moreover, $\widetilde{\mathcal{S}}$ is
invariant under the action of $M_z$. Therefore, by the
Beurling--Lax--Halmos Theorem, Theorem \ref{BLH-D}, there is an
inner function $\Theta$ in $H^\infty_{\clb(\clw)}(\mathbb{D})$ for
which $\widetilde{\mathcal{S}} = \Theta H^2_{\clw}(\mathbb{D})$.
Further, since $\widetilde{\mathcal{S}}$ has finite codimension in
$H^2_{\clw}(\mathbb{D})$ and the dimension of $\clw$ is finite, it
follows that the matrix  entries of $\Theta$ are rational functions
with poles outside the closed unit disk and $\Theta(e^{it})$ is
unitary for $e^{it}$ in $\mathbb{T}$ (cf. \cite{NF}, Chapter VI,
Section 4).

Now the determinant, $\det \Theta$, is a scalar-valued rational
inner function in $H^\infty(\mathbb{D})$ and hence is a finite
Blaschke product. Using Cramer's Rule one can show that (cf.
\cite{He}, Theorem 11) \[(\det\Theta)H^2_{\clw}(\mathbb{D})
\subseteq\Theta H^2_{\clw}(\mathbb{D}),\] which implies that
\[
\dim_{\mathbb{C}} H^2_{\clw}(\mathbb{D})/\Theta
H^2_{\clw}(\mathbb{D}) \le \dim_{\mathbb{C}}
H^2_{\clw}(\mathbb{D})/(\det \Theta)H^2_{\clw}(\mathbb{D}).
\]
Continuing, we have
\begin{align*}
\Psi(I^2_{\w_0}\cdot \mathcal{H}) &= \Psi\left(\bigvee^n_{i,j=1}
M_{z_i} M_{z_j} \mathcal{H}\right) = \bigvee^n_{i,j=1} M_{\varphi_i}
M_{\varphi_j} H^2_{\clw}(\mathbb{D})= \bigvee^n_{i=1}
M_{\varphi_i}(\Theta H^2_{\clw}(\mathbb{D}))\\ & \supseteq
\bigvee^n_{i=1} M_{\varphi_i}(\det \Theta)H^2_{\clw}(\mathbb{D})
\supseteq \bigvee^n_{i=1} \det\Theta(M_{\varphi_i} H^2_{\clw}(\mathbb{D})) = (\det\Theta) \Theta H^2_{\clw}(\mathbb{D})\\
&\supseteq (\det\Theta)^2 H^2_{\clw}(\mathbb{D}).
\end{align*}
Therefore
\[
\dim(\mathcal{H}/I^2_{\w_0}\cdot \mathcal{H}) \le \dim
H^2_{\clw}(\mathbb{D})/(\det \Theta)^2 H^2_{\clw}(\mathbb{D}).
\]
Proceeding by induction, one arrives at
\[
\dim(\mathcal{H}/I^k_{\w_0}\cdot \mathcal{H}) \le \dim
H^2_{\clw}(\mathbb{D})/(\det \Theta)^k H^2_{\clw}(\mathbb{D}),
\]
for each positive integer $k$. Also
\[
h^{\w_0}_{\mathcal{H}}(k) \le \dim
H^2_{\clw}(\mathbb{D})/(\det\Theta)^k H^2_{\clw}(\mathbb{D}) = k d
\dim \clw \quad \text{for}\quad k\ge N_{\mathcal{H}},
\]where $d$ is the dimension  of
$H^2/(\det \Theta)H^2$. Hence, the degree of
$h^{\w_0}_{\mathcal{H}}$ is at most one.\qed

Combining this theorem with Theorem \ref{pro2} yields the following
result.

\begin{Theorem}\label{thm2}
If $\mathcal{R}$ is a semi-Fredholm, quasi-free Hilbert module over
$A(\Omega)$ with $\Omega\subset \mathbb{C}^n$ having a pure
isometrically isomorphic submodule of finite codimension, then
$n=1$.
\end{Theorem}

\subsection{Hilbert modules over $A(\mathbb{D})$}
By virtue of Theorem \ref{thm2}, one can immediately reduce to the
case of domains $\Omega$ in $\mathbb{C}$ if there exists a pure
isometrically isomorphic submodule of finite codimension.

The purpose of this subsection is to prove that for a quasi-free
Hilbert module $\mathcal{R}$ over $A(\mathbb{D})$, the existence of
a pure unitarily equivalent submodule of finite codimension implies
that $\mathcal{R}$ is unitarily equivalent to
$H^2_{\mathcal{E}}(\mathbb{D})$ with $\dim \mathcal{E} = \mbox{rank}
\clr$.

\begin{Theorem}\label{thm3}
Let $\clr$ be a finite rank, quasi-free Hilbert module over
$A(\mathbb{D})$ which is semi-Fredholm for $\omega$ in $\mathbb{D}$.
Assume there exists a pure module isometry $U$ such that $\dim \clr/
U \clr <~\infty$. Then $\clr$ and  $H^2_{\cle}(\mathbb{D})$ are
$A(\mathbb{D})$-module isomorphic where $\cle$ is a Hilbert space
with $\mbox{dim~} \cle$ equal to the multiplicity of $\clr$.
\end{Theorem}

\NI\textbf{Proof.} As in Proposition \ref{pro1}, without loss of
generality one can assume that $\clr \cong H^2_{\clw}(\mathbb{D})$,
where $\clw = \clr \ominus U \clr$ with $\dim \clw < \infty$ and $U$
corresponds to $M_z$. Let $M_{\varphi}$ denote the operator on
$H^2_{\mathcal{W}}(\mathbb{D})$ unitarily equivalent to  module
multiplication by $z$ on $\mathcal{R}$, where $\varphi$ is in
$H^{\infty}_{\clb(\clw)}(\mathbb{D})$ with $\| \varphi (z) \| \leq
1$ for all $z$ in $\mathbb{D}$.

\NI Since the operator $M_{\varphi}$ is defined by module
multiplication on $H^2_{\clw}(\mathbb{D})$ and the corresponding
$A(\mathbb{D})$-module has finite rank,  it is enough to show that
$\varphi$ is inner. Hence $M_\varphi$ would be a pure isometry so
that $H^2_{\clw}(\mathbb{D})$ and $H^2(\mathbb{D})$ would be
$A(\mathbb{D})$-module isomorphic.

\NI Since the range of $M_\varphi- w I$ has finite codimension in
$H^2_{\mathcal{W}}(\mathbb{D})$, it follows that the operator
$M_{\varphi} - w I$ has closed range for each $w \in \mathbb{D}$.
Now $\mbox{ker} (M_{\varphi} - w I) = \{0\}$, by Lemma 1.1 in
\cite{DS1}, implies that $M_{\varphi} - w I$ is bounded below. Then
by Lemma 2.1 in \cite{DS1}, $(L_{\varphi} - w I)$ is bounded below
on $L^2_{\cle}(\mathbb{D})$, where $L_{\varphi}$ is the Laurent
operators with symbol $\varphi$.

\NI For each $w$ in $\mathbb{D}$ and $k$ in $\mathbb{N}$, define
$$E^{w}_k = \{ e^{it}\colon\  \mbox{dist}(\sigma(\varphi(e^{it})), w) < \frac{1}{k}\},$$
where $\sigma(\varphi(e^{it}))$ denotes the spectrum of the matrix
$\varphi(e^{it})$.

\NI Then either $\sigma(\varphi(e^{it})) \subset \mathbb{T}$ a.e or
there exists a $w_0$ in $\mathbb{D}$ such that  $m(E^{w_0}_k) > 0$
for all $k \in \mathbb{N}$. In the latter case, one can find a
sequence of functions $\{f_k\}$ in $L^2_{\mathcal{E}}(\mathbb{T})$
such that $f_k$ is supported on $E^{w_0}_k$, $\|f_k(e^{it})\| = 1$
for $e^{it}$ in $E^{\omega_0}_k$ and
$$\| \varphi(e^{it}) f_k (e^{it}) - w_0
f_k (e^{it}) \| \leq \frac{1}{k}.$$ It then follows that
$$\| (L_{\varphi} - w_0) f_k \| \leq \frac{1}{k}\|f_k\|$$ for all $k$ in $\mathbb{N}$, which contradicts the fact that
$L_{\varphi} - w_0I$ is bounded below. Hence,
$\sigma(\varphi(e^{it})) \subset \mathbb{T}$, a.e.\ and hence
$\varphi(e^{it})$ is unitary a.e. Therefore, $T_{\varphi}$ is a pure
isometry and the Hilbert module $H^2_{\cle}(\mathbb{D})$ determined
by $T_\varphi$ is $A(\mathbb{D})$-module isomorphic with
$H^2_{\cle}(\mathbb{D})$.\qed

This result can not be extended to the case in which $U$ is not
pure. For example, for $\mathcal{R} = H^2(\mathbb{D}) \oplus
L^2_a(\mathbb{D})$, one could take $U = M_z\oplus I$.

\NI \textbf{Further results and comments:}

\begin{enumerate}

\item All of the material in this section is taken from \cite{DS1}.

\item For the Bergman modules over the unit ball, one
can show (cf. \cite{C-G,P,Rich}) that no proper submodule is
unitarily equivalent to the Bergman module itself. These issues are
thoroughly discussed in \cite{S-HM}.

\item In a sense, the existence of a Hilbert module with unitarily equivalent
submodules is a rare phenomenon. The following example shows that
the problem is more complicated even in the sense of
quasi-similarity.

\NI\textsf{Example:} The Hardy module $H^2(\mathbb{D}^2)$ is not
quasi-similar to the submodule $H^2(\mathbb{D}^2)_0 = \{f \in
H^2(\mathbb{D}^2): f(0) = 0\}$ of $H^2(\mathbb{D}^2)$. Suppose $X$
and $Y$ define a quasi-affinity between $H^2(\mathbb{D}^2)$ and
$H^2(\mathbb{D}^2)_0$. Then the localized maps $X_{0}$ and $Y_{0}$
are isomorphisms between $\mathbb{C}_{0}$ and $ \mathbb{C}_0 \oplus
\mathbb{C}_0$ (see Section \ref{CHM}) which is impossible.

\item Theorem \ref{thm3} can be extended to the case of a
finitely-connected domain $\Omega$ with a nice boundary, that is,
$\Omega$ for which $\partial\Omega$ is the finite union of simple
closed curves. Here it is convenient to recall the notion of the
\textit{bundle shift} $H^2_\alpha(\Omega)$ for $\Omega$ determined
by the unitary representation $\alpha$ of the fundamental group
$\pi_1(\Omega)$ of $\Omega$. The bundle shift $H^2_\alpha(\Omega)$
is the Hardy space of holomorphic sections of the flat unitary
bundle over $\Omega$ determined by $\alpha$ (see \cite{A-D},
\cite{B-JOT}). The reader is referred to (\cite{DS1}, Theorem 2.8)
for a proof of the following theorem.

\begin{Theorem}\label{thm4}
Let $\mathcal{R}$ be a finite rank, quasi-free Hilbert module over
$A(\Omega)$, where $\Omega$ is a finitely-connected domain in
$\mathbb{C}$ with  nice boundary, which is semi-Fredholm for
$\omega$ in $\Omega$. Let $U$ be a pure module isometry such that
$\dim \mathcal{R}/U\mathcal{R} < \infty$. Then there is a unitary
representation $\alpha$ of $\pi^1(\Omega)$ on some finite
dimensional Hilbert space such that $\mathcal{R}$ and the bundle
shift $H^2_\alpha(\Omega)$, are $A(\Omega)$-module isomorphic.
\end{Theorem}

\item In \cite{DNYJ}, Douglas proved the following
result on rank one quasi-free Hilbert modules.

\begin{Theorem}Let $\clr$ be a rank one quasi-free Hilbert module
over $A(\Omega)$, where $\Omega = \mathbb{B}^n$ or $\mathbb{D}^n$.
Suppose each submodule $\cls$ of $\clr$ is isometrically isomorphic
to $\clr$. Then $n = 1$ and $\clr \cong H^2(\mathbb{D})$ and the
module map $M$ on $\clr$ is the Toeplitz operator $M_{\varphi}$,
where $\varphi$ is a conformal self map of $\mathbb{D}$ onto itself.
\end{Theorem}

\item The notion of Hilbert–-Samuel polynomials for Hilbert modules is a
relatively new concept and was introduced by Douglas and Yan in 1993
\cite{D-Y}. Because of its strong interaction with commutative
algebra and complex analytic geometry, Hilbert module approach to
Hilbert-Samuel polynomial and Samuel multiplicity has had a
spectacular development since its origin. The reader is referred to
the recent work by Eschmeier \cite{E1}, \cite{E2}, \cite{E3},
\cite{E-IJM} and Fang \cite{F1}, \cite{F2}, \cite{F3}.

\end{enumerate}

\section{Corona condition and Fredholm Hilbert modules}\label{SFHM}

The purpose of this section is to apply techniques from Taylor's
theory, in terms of Koszul complex, Berezin transforms and
reproducing kernel method to quasi-free Hilbert modules and obtain a
connection between Fredholm theory and corona condition.

\subsection{Koszul complex and Taylor invertibility}
In this subsection, the notion of Taylor's invertibility (see
\cite{Tacta}, \cite{TJFA}) for commuting tuples of operators on
Hilbert spaces will be discussed.

Let $\cle^n$ be the exterior algebra generated by $n$ symbols
$\{e_1, \ldots, e_n\}$ along with identity $e_0$, that is, $\cle^n$
is the algebra of forms in $\{e_1, \ldots, e_n\}$ with complex
coefficients and $e_i \wedge e_j = - e_j \wedge e_i$ for all $1 \leq
i, j \leq n$. Let $\cle^n_k$ be the vector subspace of $\cle^n$
generated by the basis \[\{e_{i_1} \wedge \cdots \wedge e_{i_k} :
1\leq  i_1 < \ldots < i_k \leq n\}.\] In particular, \[\cle^n_i
\wedge \cle^n_j \subseteq \cle^n_{i+j},\] and \[\cle^n = \mathbb{C}
e_1 \wedge \cdots \wedge e_n.\] Moreover \[\mbox{dim~} \cle^n_k = {n
\choose k},\] that is, $\cle^n_k$ is isomorphic to $\mathbb{C}^{n
\choose k}$ as a vector space over $\mathbb{C}$. Also note that
$\cle^n$ is graded:
\[\cle^n = \mathop{\sum}_{k= 0}^{\infty} \cle^n_k.\]Define the
\textit{creation operator} $E_i : \cle^n \raro \cle^n$, for each $1
\leq i \leq n$, by $E_i \eta = e_i \wedge \eta$ and $E_0 \eta =
\eta$ for all $\eta \in \cle^n$. In particular, note that $\cle^n$
is a finite dimensional vector space. Then the anticommutation
relation follows easily:
\[E_i E_j = - E_j E_i \quad \quad \mbox{and} \quad \quad E_i^* E_j +
E_j E^*_i = \delta_{ij} E_0.\]

Now let $T = (T_1, \ldots, T_n)$ be a commuting tuple of operators
on $\clh$. Let $\cle^n(T) = \clh \otimes_{\mathbb{C}} \cle^n$ and
$\cle^n_k(T) = \clh \otimes_{\mathbb{C}} \cle^n_k \subset \cle^n(T)$
and define $\partial_{T} \in \clb(\cle^n(T))$ by
\[\partial_{T} = \mathop{\sum}_{i=1}^n T_i \otimes E_i.\]It
follows easily from the anticommutation relationship that
$\partial_{T}^2 = 0$. The \textit{Koszul complex} $K(T)$ is now
defined to be the (chain) complex \[K(T): 0 \longrightarrow
\cle^n_0(T) \stackrel{\partial_{1, \clh}}\longrightarrow
\cle^n_{1}(T) \stackrel{\partial_{2, \clh}}\longrightarrow \cdots
\stackrel{\partial_{n-1, \clh}}\longrightarrow \cle^n_{n-1}(T)
\stackrel{\partial_{n, \clh}}\longrightarrow \cle^n_n(T)
\longrightarrow 0,\]where $\cle^n_k(T)$ is the collection of all
$k$-forms in $\cle^n(T)$ and $\partial_{k, T}$, \textit{the
differential}, is defined by \[\partial_{k, T} =
\partial_{T}|_{\cle^n_{k-1}(T)} : \cle^n_{k-1}(T) \raro \cle^n_{k}(T). \quad \quad \quad (k = 1, \ldots, n)\]
For each $k = 0, \ldots, n$ the cohomology vector space associated
to the Koszul complex $K(T)$ at $k$-th stage is the vector space
\[H^k(T) = \mbox{ker~} \partial_{k+1, T}/ \mbox{ran}
\partial_{k, T}.\]Here $\partial_{0, T}$ and $\partial_{n+1, T}$ are the zero map.
A commuting tuple of operators $T$ on $\clh$ is said to be
\textit{invertible} if $K(T)$ is exact. The \textit{Taylor spectrum}
of $T$ is defined as \[\sigma_{Tay}(T) = \{\w \in \mathbb{C}^n : K(T
- \w I_{\clh}) \mbox{~is not exact~}\}.\]The tuple $T$ is said to be
a \textit{Fredholm tuple} if
\[\mbox{dim~}\Big[H^k(T)\Big] < \infty, \quad \quad(k = 0, 1, \ldots,
n)\]and \textit{semi-Fredholm tuple} if the last cohomology group,
\[H^n(T) = \clh/\mathop{\sum}_{i=1}^n T_i \clh, \] of its Koszul
complex in finite dimensional.

\NI If $T$ is a Fredholm tuple, then the index of $T$ is
\[\mbox{ind} T:= \mathop{\sum}_{k=0}^n (-1)^k \mbox{dim~} \Big[
H^k(T)\Big].\]The tuple $T$ is said to be Fredholm (or
semi-Fredholm) at $\w \in \mathbb{C}^n$ if the tuple $T - \w
I_{\clh}$ is Fredholm (or semi-Fredholm).

Viewing the tuple $T$ as a Hilbert module over $\mathbb{C}[\z]$, it
follows that $T$ is semi-Fredholm at $\w$ if and only if
\[\dim [\clh/{I_{\w}\cdot \clh}] < \infty.\] In particular, note that
$\mathcal{H}$ semi-Fredholm at $\w$ implies that $I_{\w} \clh$ is a
closed submodule of $\mathcal{H}$.

\subsection{Weak corona property}

Let $\{\varphi_1, \ldots, \varphi_k\} \subseteq
H^{\infty}_{\clb(\cle)}(\mathbb{B}^n)$ be a $k$-tuple of commuting
$\clb(\cle)$-valued functions where $\cle$ is a Hilbert space. Then
the tuple is said to have the \textit{weak corona property} if there
exists an $\epsilon
>0$ and $1>\delta > 0$ such that
$$\sum_{i=1}^{k} \varphi_i(\z) \varphi_i(\z)^* \geq \epsilon I_{\cle},$$for all
$\z$ satisfying $1>\|\z\| \geq 1-\delta$.

\NI The tuple $\{\varphi_1, \ldots, \varphi_k\}$  is said to have
the \textit{corona property} if \[\sum_{i=1}^k \varphi_i(\z)
\varphi_i(\z)^* \geq \epsilon I_{\cle},\] for all $\z \in
\mathbb{B}^n$.

For n=l and $\cle = \mathbb{C}$, the Carleson's corona theorem (see
\cite{CAM}) asserts that:

\begin{Theorem}\textsf{(Carleson)} A set $\{\varphi_1, \ldots, \varphi_k\}$ in
$H^{\infty}(\mathbb{D})$ satisfies $\sum_{i=1}^{k} |\varphi_i(z)|
\geq \epsilon$ for all $z$ in $\mathbb{D}$  for some $\epsilon >0$
if and only if there exist $\{\psi_1, \ldots, \psi_k\} \subset
H^{\infty}(\mathbb{D})$ such that \[\sum_{i=1}^{k} \varphi_i \psi_i
= 1.\]
\end{Theorem}

Also one has the following fundamental result of Taylor (see
\cite{TJFA}, Lemma 1):

\begin{Lemma}\label{Tay}
Let $(T_1, \ldots, T_k)$ be in the center of an algebra $\cla$
contained in $\cll(\clh)$ such that there exists $(S_1, \ldots,
S_k)$ in $\cla$  satisfying $\sum_{i=1}^{k} T_i S_i = I_{\clh}$.
Then the Koszul complex for $(T_1, \ldots, T_k)$ is exact.
\end{Lemma}

Now consider a contractive quasi-free Hilbert module $\clr$ over
$A(\mathbb{D})$ of multiplicity one, which therefore has
$H^{\infty}(\mathbb{D})$ as the multiplier algebra.

\begin{Proposition}
Let $\clr$ be a contractive quasi-free Hilbert module over
$A(\mathbb{D})$ of multiplicity one and $\{\varphi_1, \ldots,
\varphi_k\}$ be a subset of $H^{\infty}(\mathbb{D})$. Then the
Koszul complex for the $k$-tuple $(M_{\varphi_1}, \ldots,
M_{\varphi_k})$ on $\clr$ is exact if and only if $\{\varphi_1,
\ldots, \varphi_k\}$ satisfies the corona property.
\end{Proposition}
\NI \textsf{Proof.} If $\sum_{i=1}^{k} \varphi_i \psi_i = 1$ for
some $\{\psi_1, \ldots, \psi_k\} \subset H^{\infty}(\mathbb{D})$,
then the fact that $M_{\Phi}$ is Taylor invertible follows from
Lemma \ref{Tay}. On the other hand, the last group of the Koszul
complex is $\{0\}$ if and only if the row operator $M_{\varphi}$ in
$\clb(\clr^k, \clr)$ is bounded below which, as before, shows that
$\sum_{i=1}^{k} |\varphi_i(z)|$ is bounded below on $\mathbb{D}$.
This completes the proof. \qed

The missing step to extend the result from $\mathbb{D}$ to the open
unit ball $\mathbb{B}^n$ is the fact that it is unknown if the
corona condition for $\{\varphi_1, \ldots, \varphi_k\}$ in
$H^{\infty}(\mathbb{B}^n)$ is equivalent to the Corona property.

\subsection{Semi-Fredholm implies weak corona}

Let $\clh_K$ be a scalar-valued reproducing kernel Hilbert space
over $\Omega$ and $F \in \clb(\clh_K)$. Then the \textit{Berezin
transform} (see \cite{DD}) of $F$ is denoted by $\hat{F}$ and
defined by
\[\hat{F}(\z) = \langle F \frac{K(\cdot, \z)}{\|K(\cdot, \z)\|},
\frac{K(\cdot, \z)}{\|K(\cdot, \z)\|}\rangle. \quad\quad (\z \in
\Omega)\]

Note that the multiplier space of a rank one quasi-free Hilbert
module $\clr$ over $A(\mathbb{B}^n)$ is precisely
$H^{\infty}(\mathbb{B}^n)$, since $\clr$ is the completion of
$A(\mathbb{B}^n)$, by definition (see Proposition 5.2 in \cite{DD}).

\begin{Theorem}\label{prop1}
Let $\clr$ be a contractive quasi-free Hilbert module over
$A(\mathbb{B}^n)$ of multiplicity one and $\{\varphi_1, \ldots,
\varphi_k\}$ be a subset of $H^{\infty}(\mathbb{B}^n)$. If
$(M_{\varphi_1}, \ldots, M_{\varphi_k})$ is a semi-Fredholm tuple,
then $\{\varphi_1, \ldots, \varphi_k\}$ satisfies the weak corona
condition.
\end{Theorem}
\NI\textsf{Proof.} Let $K : \mathbb{B}^n \times \mathbb{B}^n \raro
\mathbb{C}$ be the kernel function for the quasi-free Hilbert module
$\clr$. By the assumption, the range of the row operator $M_{\Phi} =
(M_{\varphi_1}, \ldots, M_{\varphi_k}): \clr^k \raro \clr$ in $\clr$
has finite co-dimension, that is, $$\text{dim} [\clr/(M_{\varphi_1}
\clr + \ldots + M_{\varphi_k} \clr)] < \infty,$$ and, in particular,
$M_{\Phi}$ has closed range. Consequently, there is a finite rank
projection $F$ in $\clb(\clr)$ such that
$$M_{\Phi} M_{\Phi}^* + F = \sum_{i=1}^{k} M_{\varphi_i}
M_{\varphi_i}^* + F : \clr \raro \clr$$ is bounded below. Therefore,
there exists a $c>0$ such that $$\langle F K(\cdot, \z), K(\cdot,
\z)\rangle +  \langle \sum_{i=1}^{k} M_{\varphi_i} M_{\varphi_i}^*
K(\cdot, \z), K(\cdot, \z)\rangle \; \geq c \|K(\cdot, \z)\|^2,$$
for all $\z \in \mathbb{B}^n$. Therefore,
$$\|K(\cdot, \z)\|^2 \hat{F}(\z) + \|K(\cdot, \z)\|^2 \Big(\sum_{i=1}^{k}
\varphi_i(\z) \varphi_i^* (\z)\Big) \geq c \|K(\cdot, \z)\|^2,$$ and
so
$$ \hat{F}(\z) + \sum_{i=1}^{k} \varphi_i(\pmb{z})
\varphi_i(\pmb{z})^* \geq c,$$ for all $\pmb{z}$ in $\mathbb{B}^n$.
Using the known boundary behavior of the Berezin transform (see
Theorem 3.2 in \cite{DD}), since $F$ is finite rank we have that
$|\hat{F}(\z)| \leq \frac{c}{2}$ for all $\z$ such that $1 > \|\z\|
> 1- \delta$ for some $1 > \delta>0$ depending on $c$. Hence
$$\sum_{i=1}^{k} \varphi_i(\z) \varphi_i(\z)^* \geq
\frac{c}{2},$$ for all $\z$ such that $1>\|\z\|>1-\delta
>0$, which completes the proof. \qed

The key step in this proof is the vanishing of the Berezin transform
at the boundary of $\mathbb{B}^n$ for a compact operator. The proof
of this statement depends on the fact that $\frac{K(\cdot,
\z)}{\|K(\cdot, \z)\|}$ converges weakly to zero as $\z$ approaches
the boundary which rests on the fact that $\clr$ is contractive.

\subsection{A Sufficient condition}

\begin{Theorem}\label{THM1}
Let $\clr$ be a contractive quasi-free Hilbert module over
$A(\mathbb{D})$ of multiplicity one, which is semi-Fredholm at each
point $z$ in $\mathbb{D}$. If $\{\varphi_1, \ldots, \varphi_k\}$ is
a subset of $H^{\infty}(\mathbb{D})$, then the $k$-tuple $M_{\Phi} =
(M_{\varphi_1}, \ldots, M_{\varphi_k})$ is semi-Fredholm if and only
if it is Fredholm if and only if  $(\varphi_1, \ldots, \varphi_k)$
satisfies the weak corona condition.
\end{Theorem}
\NI \textsf{Proof.} If $M_{\Phi}$ is semi-Fredholm, then by
Proposition \ref{prop1} there exist $\epsilon >0$ and $1 > \delta>0$
such that $$\sum_{i=1}^{k} |\varphi_i(z)|^2 \geq \epsilon,$$ for all
$z$ such that $1>|z|>1-\delta > 0$. Let $\clz$ be the set $$\clz =
\{ z \in \mathbb{D} : \varphi_i (z) = 0 \; \text{for all}\; i = 1,
\ldots, k\}.$$ Since the functions $\{\varphi_i\}_{i=1}^{k}$ can not
all vanish for $z$ satisfying $1>|z|>1-\delta$, it follows that the
cardinality of the set $\clz := N$ is finite. Let
$$\clz = \{ z_1, z_2, \ldots, z_N\}$$ and $l_j$ be the smallest
order of the zero at $z_j$ for all $\varphi_j$ and $1 \leq j \leq
k$. Let $B(z)$ be the finite Blaschke product with zero set
precisely $\clz$ counting the multiplicities. Note that $\xi_i : =
\frac{\varphi_i}{B} \in H^{\infty}(\mathbb{D})$ for all $i=1,
\ldots, k$. Since $\{\varphi_1, \ldots, \varphi_k\}$ satisfies the
weak corona property, it follows that $\sum_{i=1}^{k} |\xi_i (z)|^2
\geq \epsilon$ for all $z$ such that $1>|z|>1-\delta$. Note that
$\{\xi_1, \ldots, \xi_n\}$ does not have any common zero and so
$\sum_{i=1}^{k} |\xi_i (z)|^2 \geq \epsilon$, for all $z$ in
$\mathbb{D}$. Therefore, $\{\xi_1, \ldots, \xi_k\}$ satisfies the
corona property and hence there exists $\{\eta_1, \ldots, \eta_k\}$,
a subset of $H^{\infty}(\mathbb{D})$, such that \[\sum_{i=1}^{k}
\xi_i(z) \eta_i(z) = 1,\] for all $z$ in $\mathbb{D}$. Thus,
\[\sum_{i=1}^{k} \varphi_i(z) \eta_i(z) = B,\] for all $z$ in
$\mathbb{D}$. This implies $\sum_{i=1}^{k} M_{\varphi_i} M_{\eta_i}
= M_B$, and consequently, $\sum_{i=1}^{k} \overline{M}_{\varphi_i}
\overline{M}_{\eta_i} = \overline{M}_B,$ where
$\overline{M}_{\varphi_i}$ is the image of $M_{\varphi_i}$ in the
Calkin algebra, $\clq(\clr) = \clb(\clr)/\clk(\clr)$. But the
assumption that $M_{z-w}$ is Fredholm for all $w$ in $\mathbb{D}$
yields that $M_B$ is Fredholm. Therefore,  $X = \sum_{i=1}^{k}
\overline{M}_{\varphi_i} \overline{M}_{\eta_i}$ is invertible.
Moreover, since $X$ commutes with the set
\[\{\overline{M}_{\varphi_1}, \ldots, \overline{M}_{\varphi_k},
\overline{M}_{\eta_1}, \ldots, \overline{M}_{\eta_k}\},\] it follows
that $(M_{\varphi_1}, \ldots, M_{\varphi_k})$ is a Fredholm tuple,
which completes the proof. \qed

\vspace{0.2in}

\NI \textbf{Further results and comments:}

\begin{enumerate}
\item In Theorem 8.2.6 in \cite{EP}, a version of Theorem \ref{prop1} is
established in case $\clr$ is the Bergman module on $\mathbb{B}^n$.

\item The converse of Theorem \ref{prop1} is known for the Bergman space
for certain domains in $\mathbb{C}^n$ (see Theorem 8.2.4 in
\cite{EP} and pages 241-242). A necessary condition for the converse
to hold for the situation in Theorem \ref{prop1} is for the
$n$-tuple of co-ordinate multiplication operators to have essential
spectrum equal to $\partial \mathbb{B}^n$, which is not automatic,
but is true for the classical spaces.

\item One prime reason to establish a converse, in Theorem
\ref{THM1}, is that one can represent the zero variety of the ideal
generated by the functions in terms of a single function, the finite
Blaschke product (or polynomial). This is not surprising since
$\mathbb{C}[z]$ is a principal ideal domain.

\item As pointed out in the monograph by Eschmeier and Putinar,
the relation between corona problem and the Taylor spectrum is not
new (cf. \cite{LH}, \cite{RW}).

\item This section is mainly based on \cite{DS3} and closely related to
\cite{DE12} and \cite{DS2}.

\item In \cite{Ven}, Venugopalkrishna developed a Fredholm theory
and index theory for the Hardy module over strongly pseudoconvex
domains in $\mathbb{C}^n$.

\item An excellent source of information concerning Taylor spectrum is the monograph by Muller
\cite{M-book}. See also the paper \cite{Curto} and the survey
\cite{C-S} by Curto and the book by Eschmeier and Putinar \cite{EP}.

\end{enumerate}

\section{Co-spherically contractive Hilbert modules}\label{CSCHM}

A Hilbert module over $\mathbb{C}[\bm{z}]$ is said to be {\it
co-spherically contractive}, or define a  {\it row contraction}, if
$$\|\sum_{i=1}^{n} M_i h_i\|^2 \leq \sum_{i=1}^{n} \|h_i\|^2, \quad
(h_1, \ldots, h_n \in \clh),$$ or, equivalently, if $\sum_{i=1}^{n}
M_i M_i^* \leq I_{\clh}$. Define the defect operator and the defect
space of $\clh$ as \[D_{*\clh} = (I_{\clh} - \sum_{i=1}^n M_i
M_i^*)^{\frac{1}{2}} \in \cll(\clh),\] and \[\cld_{*\clh} =
\overline{\mbox{ran}} D_{*\clh},\] respectively. We denote
$D_{*\clh}$ and $\cld_{*\clh}$ by $D_*$ and $\cld_*$ respectively,
if $\clh$ is clear from the context.

If $n = 1$ then $\clh$ is a contractive Hilbert module over
$A(\mathbb{D})$ (see Section \ref{CHM}).

\subsection{Drury-Arveson Module}
Natural examples of co-spherically contractive Hilbert modules over
$\mathbb{C}[\bm{z}]$ are the Drury-Arveson module, the Hardy module
and the Bergman module, all defined on $\mathbb{B}^n$.

One can identify the Hilbert tensor product $H^2_n \otimes \cle$
with the $\cle$-valued $H^2_n$ space $H^2_n(\cle)$ or the
$\clb(\cle)$-valued reproducing kernel Hilbert space with kernel
function $(\z, \w) \mapsto (1 - \sum\limits_{i=1}^{n} z_i
\bar{w}_i)^{-1} I_{\cle}$. Then
$$H^2_n(\cle) = \{ f \in \clo (\mathbb{B}^n, \cle): f(z) =
\sum_{\bm{k} \in \mathbb{N}^n} a_{\bm{k}} z^{\bm{k}}, a_{\bm{k}} \in
\cle, \|f\|^2 : = \sum_{\bm{k} \in \mathbb{N}^n} \frac{ \|
a_{\bm{k}} \|^2}{\gamma_{\bm{k}}} < \infty \},$$ where
$\gamma_{\bm{k}} = \frac{(k_1 + \cdots + k_n)!}{k_1 ! \cdots k_n!}$
are the multinomial coefficients and $\bm{k} \in \mathbb{N}^n$.

Given a co-spherically contractive Hilbert module $\clh$, define the
completely positive map $P_{\clh} : \cll(\clh) \rightarrow
\cll(\clh)$ by \[P_{\clh} (X) = \sum_{i=1}^{n} M_i X M^*_i,\] for
all $X \in \cll(\clh)$. Note that
\[I_{\clh} \geq P_{\clh} (I_{\clh}) \geq P^2_{\clh} (I_{\clh}) \geq
\cdots \geq P^l_{\clh} (I_{\clh}) \geq \cdots \geq 0.\] In
particular,
$$P_{\infty}(\clh) := \mbox{SOT} - \mathop{\lim}_{l \raro \infty}
P_{\clh}^l (I_{\clh})$$ exists and $0 \leq P_{\infty}(\clh) \leq
I_{\clh}$. The Hilbert module $\clh$ is said to be \textit{pure} if
$$P_{\infty}(\clh) = 0.$$  Examples of pure co-spherically
contractive Hilbert modules over $\mathbb{C}[\bm{z}]$ includes the
submodules and quotient modules of vector-valued Drury-Arveson
module.

\subsection{Quotient modules of $H^2_n(\cle)$}

First recall a standard result from algebra: Any module is
isomorphic to a quotient of a free module. The purpose of this
subsection is to prove an analogous result for co-spherically
contractive Hilbert modules: any pure co-spherically contractive
Hilbert module is isomorphic to a quotient module of the
Drury-Arveson module with some multiplicity.

\begin{Theorem}\label{PIH}
Let $\clh$ be a co-spherically contractive Hilbert module over
$\mathbb{C}[\bm{z}]$. Then there exists a unique co-module map
$\Pi_{\clh} : \clh \raro H^2_n(\cld_*)$ such that
\[(\Pi_{\clh} h)(\bm{w}) = D_* (I_{\clh} - \mathop{\sum}_{i=1}^n w_i
M_i^*)^{-1}h, \quad \quad(\bm{w} \in \mathbb{B}^n,\,h \in \clh)\]
and $\Pi_{\clh}^* \Pi_{\clh} = I_{\clh} - P_{\infty}(\clh)$.
Moreover, $\Pi_{\clh}^* (k_n(\cdot, \bm{w}) \eta) =  (I_{\clh} -
\sum_{i=1}^n \bar{w}_i M_i)^{-1} D_* \eta$ for $\w \in \mathbb{B}^n$
and $\eta \in \cld_*$.
\end{Theorem}

\NI\textsf{Proof.} First, note that for each $\bm{w} \in
\mathbb{B}^n$ that \[\begin{split}\|\mathop{\sum}_{i=1}^n w_i
M_i^*\| & = \|(w_1 I_{\clh}, \ldots, w_n I_{\clh})^* (M_1, \ldots,
M_n)\|\leq \|(w_1 I_{\clh}, \ldots, w_n I_{\clh})^*\| \|(M_1,
\ldots, M_n)\| \\ & = (\sum_{i=1}^n |w_i|^2)^{\frac{1}{2}}
\|\sum_{i=1}^n M_i M_i^*\|^{\frac{1}{2}} = \|\bm{w}\|_{\mathbb{C}^n}
\|\sum_{i=1}^n M_i M_i^*\|^{\frac{1}{2}} < 1.\end{split}\]Therefore,
$\Pi_{\clh} : \clh \raro H^2_n(\cld_*)$ defined by
\[(\Pi_{\clh} h)(\bm{z}):= D_* (I_{\clh} - \mathop{\sum}_{i=1}^n z_i
M_i^*)^{-1} h = \sum_{\bm{k} \in \mathbb{N}^n} (\gamma_{\bm{k}} D_*
M^{*\bm{k}} h) z^{\bm{k}},\]for all $h \in \clh$ and $\bm{z} \in
\mathbb{B}^n$, is a bounded linear map. Also the equalities
\[\begin{split} \|\Pi_{\clh} h\|^2 & = \|\sum_{\bm{k} \in
\mathbb{N}^n} (\gamma_{\bm{k}} D_* M^{*\bm{k}} h) z^{\bm{k}}\|^2 =
\sum_{\bm{k} \in \mathbb{N}^n} \gamma_{\bm{k}}^2 \|D_* M^{*\bm{k}}
h\|^2 \|z^{\bm{k}}\|^2 = \sum_{\bm{k} \in
\mathbb{N}^n}\gamma_{\bm{k}}^2 \|D_* M^{*\bm{k}} h\|^2
\frac{1}{\gamma_{\bm{k}}}\\ & = \sum_{\bm{k} \in
\mathbb{N}^n}\gamma_{\bm{k}} \|D_* M^{*\bm{k}} h\|^2 =
\sum_{l=0}^\infty \sum_{|\bm{k}| = l}\gamma_{\bm{k}} \|D_*
M^{*\bm{k}} h\|^2 = \sum_{l=0}^\infty
\sum_{|\bm{k}|=l}\gamma_{\bm{k}} \langle M^{\bm{k}} D^2_*
M^{*\bm{k}} h, h \rangle
\\ & = \sum_{l=0}^\infty \langle \sum_{|\bm{k}|=l}\gamma_{\bm{k}}  M^{\bm{k}} D^2_* M^{*\bm{k}} h, h
\rangle = \sum_{l=0}^\infty \langle P^l_{\clh}(D^2_*) h, h \rangle =
\sum_{l=0}^\infty \langle P^l_{\clh}(I_{\clh} - P_{\clh}(I_{\clh}))
h, h \rangle \\& = \sum_{l=0}^\infty \langle (P^l_{\clh}(I_{\clh}) -
P^{l+1}_{\clh}(I_{\clh})) h, h \rangle = \sum_{l=0}^\infty (\langle
P^l_{\clh}(I_{\clh}) h, h \rangle - \langle P^{l+1}_{\clh}(I_{\clh})
h, h \rangle)\\ & = \|h\|^2 - \langle P_\infty(\clh) h,
h\rangle,\end{split}\]holds for all $h \in \clh$, where the last
equality follows from the fact that
$\{P^l_{\clh}(I_{\clh})\}_{l=0}^\infty$ is a decreasing sequence of
positive operators and that $P^0_{\clh}(I_{\clh}) = I_{\clh}$ and
$P_\infty(\clh) = \lim_{l \raro \infty} P^l_{\clh}(I_{\clh})$.
Therefore, $\Pi_{\clh}$ is a bounded linear operator and
\[\Pi_{\clh}^* \Pi_{\clh} = I_{\clh} - P_\infty(\clh).\]On the other
hand, for all $h \in \clh$ and $\bm{w} \in \mathbb{B}^n$ and $\eta
\in \cld_*$, it follows that \[\begin{split}\langle \Pi_\clh^*
(k_n(\cdot, \bm{w}) \eta), h \rangle_\clh & = \langle k_n(\cdot,
\bm{w}) \eta, D_* (I_{\clh} - \mathop{\sum}_{i=1}^n w_i M_i^*)^{-1}
h\rangle_{H^2_n(\cld_*)}\\& = \langle \sum_{\bm{k} \in \mathbb{N}^n}
(\gamma_{\bm{k}} {\bar{w}}^{\bm{k}} \eta) z^{\bm{k}}, \sum_{\bm{k}
\in \mathbb{N}^n} (\gamma_{\bm{k}} D_* M^{*\bm{k}} h) z^{\bm{k}}
\rangle_{H^2_n(\cld_*)} = \sum_{\bm{k} \in \mathbb{N}^n}
\gamma_{\bm{k}} \bar{w}^{\bm{k}} \langle M^{\bm{k}} D_* \eta, h
\rangle_{\clh}\\ & = \langle (I_{\clh} - \sum_{i=1}^n \bar{w}_i
M_i)^{-1} D_* \eta, h\rangle_{\clh},\end{split}\]that is,
\[\Pi_{\clh}^* (k_n(\cdot, \bm{w}) \eta) =  (I_{\clh} - \sum_{i=1}^n
\bar{w}_i M_i)^{-1} D_* \eta.\]Also for all $\eta \in \cld_*$ and
$\bm{l} \in \mathbb{N}^n$, \[\langle \Pi_{\clh}^* (z^{\bm{l}} \eta),
h\rangle = \langle z^{\bm{l}} \eta, \sum_{\bm{k} \in \mathbb{N}^n}
(\gamma_{\bm{k}} D_* M^{*\bm{k}} h) z^{\bm{k}}\rangle =
\gamma_{\bm{l}} \|z^{\bm{l}}\|^2 \langle \eta, D_* M^{* \bm{l}} h
\rangle = \langle M^{\bm{l}} D_* \eta, h\rangle,
\]and hence $\Pi_{\clh}$ is a co-module map.
Finally, uniqueness of $\Pi_{\clh}$ follows from the fact that
$\{z^{\bm{k}} \eta: \bm{k} \in \mathbb{N}^n, \eta \in \cld_*\}$ is a
total set of $H^2_n(\cld_*)$. This completes the proof. \qed

It is an immediate consequence of this result that if $\clh$ is a
pure co-spherical contractive Hilbert module over
$\mathbb{C}[\bm{z}]$, then $P_\infty(\clh) = 0$. Equivalently, that
$\Pi_{\clh}$ is an isometry. This yields the dilation result for
pure co-spherical contractive Hilbert modules over $\mathbb{C}[\z]$.

\begin{Corollary}\label{pure-H2n}
Let $\clh$ be a pure co-spherical contractive Hilbert module over
$\mathbb{C}[\bm{z}]$. Then \[\clh \cong \clq,\]for some quotient
module $\clq$ of $H^2_n(\cld_*)$.
\end{Corollary}
\NI\textsf{Proof.} By Theorem \ref{PIH}, the co-module map
$\Pi_{\clh}: \clh \raro H^2_n(\cld_*)$ is an isometry. In
particular, $\clq = \Pi_\clh \clh$ is a quotient module of
$H^2_n(\cld_*)$. This completes the proof. \qed

A Hilbert module $\clh$ over $\mathbb{C}[\z]$ is said to be
\textit{spherical Hilbert module} if $M_i$ is normal operator for
each $1 \leq i \leq n$ and
\[\mathop{\sum}_{i=1}^n M_i M_i^* = I_{\clh}.\]Given a spherical
Hilbert module $\clh$ over $\mathbb{C}[\z]$, there exists a unique
unital $*$-representation $\pi : C^*(\partial \mathbb{B}^n) \raro
\clb(\clh)$ such that $\pi(z_i) = M_i$ and vice versa (see
\cite{A-90}, \cite{A-92}, \cite{AAM}).

The following dilation theorem is due to Arveson \cite{AAM}.

\begin{Theorem}\label{DA-dilation}
Let $\clh$ be a co-spherical contractive Hilbert module over
$\mathbb{C}[\z]$. Then there exists a spherical Hilbert module
$\cls$ over $\mathbb{C}[\z]$ such that $H^2_n(\cld_{*}) \oplus \cls$
is a dilation of $\clh$. Equivalently, there exists a spherical
Hilbert module $\cls$ over $\mathbb{C}[\z]$ and a co-module isometry
$U: \clh \raro H^2_n(\cld_{*}) \oplus \cls$. In particular,
\[\clh \cong \clq,\]for some quotient module $\clq$ of
$H^2_n(\cld_{*}) \oplus \cls$. Moreover, the minimal dilation is
unique.
\end{Theorem}

\subsection{Curvature Inequality}

The purpose of this subsection is to compare the curvatures of the
bundles $E^*_{\clq}$ associated with a quotient module $\clq = \clh
\otimes \cle/ \cls \in B^*_m(\Omega)$ and $E^*_{\clh}$, where $\clh
\in B_1^*(\Omega)$ and $\cle$, a coefficient Hilbert space. First,
we need to recall some results from complex geometry concerning
curvatures of sub-bundles and quotient bundles (cf. \cite{GH}, pp.
78-79).

Let $E$ be a Hermitian anti-holomorphic bundle over $\Omega$
(possibly infinite rank) and $F$ be an anti-holomorphic sub-bundle
of $E$ such that the quotient $Q = E/F$ is also anti-holomorphic.
Let $\triangledown_E$ denote the Chern connection on $E$ and
$\clk_E$ the corresponding curvature form. There are two canonical
connections that we can define on $F$ and the quotient bundle $Q$.
The first ones are the Chern connections $\triangledown_F$ and
$\triangledown_Q$ on $F$ and $Q$, respectively. To obtain the second
connections, let $P$ denote the projection-valued bundle map of $E$
so that $P(\z)$ is the orthogonal projection of $E(\z)$ onto
$F(\z)$. Then
$$\triangledown_{PE} = P \triangledown_E P \quad \quad
\mbox{and}\quad \quad \triangledown_{P^{\perp}E} = P^{\perp}
\triangledown_E P^{\perp},$$ define connections on $F$ and $Q$,
respectively, where $P^{\perp} = I - P$ and $Q$ is identified fiber
wise with $P^{\perp} E$. The following result from complex geometry
relates the curvatures for these pairs of connections.
\begin{Theorem}
If $F$ is an anti-holomorphic sub-bundle of the anti-holomorphic
bundle $E$ over $\Omega$ such that $E/F$ is anti-holomorphic, then
the curvature functions for the connections $\triangledown_F, \,
\triangledown_{PE}, \, \triangledown_Q$ and
$\triangledown_{P^{\perp}E}$ satisfy $$\clk_F (\w) \geq
\clk_{PE}(\w) \quad \quad \mbox{and}\quad \quad \clk_Q(\w) \leq
\clk_{P^{\perp} E}(\w). \quad \quad (\w \in \Omega).$$
\end{Theorem}

The proof is essentially a matrix calculation involving the
off-diagonal entries of $\triangledown_E$, one of which is the
second fundamental form and the other its dual (cf. \cite{GH}).
(Note in \cite{GH}, $E$ is finite rank but the proof extends to the
more general case.)

An application of this result to Hilbert modules yields the
following:

\begin{Theorem}
Let $\clh \in B_1^*(\Omega)$ be a Hilbert module over $A(\Omega)$
(or over $\mathbb{C}[\z]$) and $\cls$ be a submodule of $\clh
\otimes \cle$ for a Hilbert space $\cle$ such that the quotient
module $\mathcal Q = (\clh \otimes \cle)/\cls$ is in
$B_m^*(\Omega)$. If $E^*_{\clh}$ and $E^*_{\mathcal Q}$ are the
corresponding Hermitian anti-holomorphic bundles over $\Omega$, then
$$P^{\perp}(\w) ( \clk_{E^*_{\clh}}(\w) \otimes I_{\cle})
P^{\perp}(\w) \geq \clk_{E^*_{\mathcal Q}}(\w). \quad (\w \in
\Omega)$$
\end{Theorem}
\NI\textsf{Proof.} The result follows from the previous theorem by
setting $E = E^*_{\clh} \otimes \cle, \, F = E^*_\cls$ and $Q =
E^*_\mathcal Q$. \qed

In particular, one has the following extremal property of the
curvature functions.

\begin{Theorem}\label{curv-ext}
Let $\clh \in B^*_m(\Omega)$ be a Hilbert module over $A(\Omega)$.
If $\clh$ is dilatable to $\clr \otimes \cle$ for some Hilbert space
$\cle$, then $$\clk_{E^*_{\clr}}(\w) \otimes I_{\cle} \geq
\clk_{E^*_{\clh}}(w).\quad \quad (\w \in \Omega)$$
\end{Theorem}

The following factorization result is a special case of Arveson's
dilation result (see Corollary 2 in \cite{DMS} for a proof).

\begin{Theorem}\label{DAkernel}
Let $\clh_k$ be a reproducing kernel Hilbert module over
$\mathbb{C}[\z]$ with kernel function $k$ over $\mathbb{B}^n$. Then
$\clh_k$ is co-spherically contractive if and only if the function
$(1 - \sum_{i=1}^n z_i \bar{w}_i)k(\z,\w)$ is positive definite.
\end{Theorem}

The following statement is now an easy consequence of Theorem
\ref{DAkernel}.

\begin{Corollary}\label{DA-curv} Let $\clh_k$ be a co-spherically
contractive reproducing kernel Hilbert module over $\mathbb{B}^n$.
Then\[\clk_{E^*_{H^2_n}} - \clk_{E^*_{\clh_k}} \geq 0.\]
\end{Corollary}

\NI \textbf{Further results and comments:}

\begin{enumerate}

\item The Drury-Arveson space has been used, first in connection
with the models for commuting contractions by Lubin in 1976
\cite{L76} (see also \cite{L77}), and then by Drury in 1978 in
connection with the von Neumann inequality for commuting contractive
tuples. However, the Drury-Arveson space has been popularized by
Arveson in 1998 \cite{AAM}.

\item The proof of Theorem \ref{PIH} is a classic example of
technique introduced by Rota \cite{Ro} in the context of similarity
problem for strict contractions. In \cite{B-PAMS}, J. Ball obtained
a several-variables analogue of Rota's model. In connection with
Rota's model, see also the work by Curto and Herrero \cite{CuHe}.

\item The converse of Theorem \ref{DA-curv} is false in general. A
converse of Theorem \ref{DA-curv} is related to the notion of
infinite divisibility (see \cite{BKM}).

\item Theorem \ref{DA-curv} is from \cite{DMS}.  For $n = 1$,
this result was obtained by Misra in \cite{GM} and was further
generalized by Uchiyama in \cite{U}.

\item Theorem \ref{DA-dilation} was proved independently by many authors (see
\cite{MV}, \cite{GP99}). Most probably, the existence of dilation
was proved for the first time by Jewell and Lubin in \cite{JL} and
\cite{L76}. However, the uniqueness part of the minimal dilation is
due to Arveson.

\item The inequality in Theorem \ref{DA-curv} shows in view of Theorem
\ref{curv-ext} that the module $H^2_n$ is an extremal element in the
set of co-spherically contractive Hilbert modules over the algebra
$\mathbb{C}[\z]$. Similarly, for the polydisk $\mathbb D^n$, the
Hardy module is an extremal element in the set of those modules over
the algebra $A(\mathbb D^n)$ which admit a dilation to the Hardy
space $H^2(\mathbb D^n) \otimes \mathcal  E$.

\item We refer the reader to Athavale \cite{A-92}, \cite{A-90} for an analytic approach and Attele and Lubin \cite{AL}
for a geometric approach to the (regular unitary) dilation theory.
In particular, Athavale proved that a spherical isometry must be
subnormal. Other related work concerning dilation of commuting
tuples of operators appears in \cite{RS}, \cite{CV1}, \cite{CV2}.

\item Motivated by the Gauss-Bonnet theorem and the curvature
of a Riemannian manifold, in \cite{AJRA} Arveson introduced a notion
of curvature which is a numerical invariant. His notion of curvature
is related to the Samuel multiplicity \cite{F-05}, Euler
characteristic \cite{AJRA} and Fredholm index \cite{GRS}.

\end{enumerate}

\end{document}